\documentclass[ejs,preprint]{imsart}

\usepackage{amsfonts}
\usepackage{mathtools} 
\usepackage{amssymb, amsmath}
\usepackage{multicol}
\usepackage{multirow}
\usepackage[mathscr]{eucal}
\usepackage{graphics}
\usepackage{graphicx}
\usepackage{verbatim}
\usepackage{color}
\usepackage{float}
\usepackage{url}

\usepackage{color}

\usepackage[noadjust]{cite}

    \newcommand{\argmin}{\mathop{\rm argmin}}

    \newcommand{\wh}{\widehat}

    \newcommand{\wt}{\widetilde}

    \newtheorem{theorem}{Theorem}[section]
    \newtheorem{lemma}{Lemma}[section]
    
    \newtheorem{proposition}{Proposition}[section]
    
    \newtheorem{corollary}{Corollary}[section]

    \newtheorem{remark}{Remark}[section]
    \numberwithin{equation}{section}

    \newenvironment{proof}[1][Proof]{\begin{trivlist}
    \item[\hskip \labelsep {\bfseries #1}]}{\end{trivlist}}

\allowdisplaybreaks

\makeatletter
\renewcommand{\@biblabel}[1]{}
\renewenvironment{thebibliography}[1]
     {\section*{\refname}%
      \@mkboth{\MakeUppercase\refname}{\MakeUppercase\refname}%
      \list{}%
           {\labelwidth=0pt
            \labelsep=0pt
            \leftmargin1.5em
            \itemindent=-1.5em
            \advance\leftmargin\labelsep
            \@openbib@code
            }%
      \sloppy
      \clubpenalty4000
      \@clubpenalty \clubpenalty
      \widowpenalty4000%
      \sfcode`\.\@m}
\makeatother
\usepackage{breakcites}

\doi{10.1214/154957804100000000}
\pubyear{0000}
\volume{0}
\firstpage{0}
\lastpage{0}

\startlocaldefs
\endlocaldefs

\begin{document}

\begin{frontmatter}

\title{Asymptotics and Optimal Bandwidth for Nonparametric Estimation of Density Level Sets}
\runtitle{Optimal Bandwidth for Level Sets}

\begin{aug}
  \author{Wanli Qiao\thanksref{t2}\ead[label=e1]{wqiao@gmu.edu}}

  \address{ Department of Statistics\\
George Mason University\\
4400 University Drive, MS 4A7\\
Fairfax, VA 22030\\
           \printead{e1}}

  \thankstext{t2}{Partially supported by NSF grants DMS 1821154 and FET 1900061, and a Jeffress Memorial Trust Award.}
  \runauthor{W. Qiao}

\end{aug}

\begin{abstract}
Bandwidth selection is crucial in the kernel estimation of density level sets. A risk based on the symmetric difference between the estimated and true level sets is usually used to measure their proximity. In this paper we provide an asymptotic $L^p$ approximation to this risk, where $p$ is characterized by the weight function in the risk. In particular the excess risk corresponds to an $L^2$ type of risk, and is adopted to derive an optimal bandwidth for nonparametric level set estimation of $d$-dimensional density functions ($d\geq 1$). A direct plug-in bandwidth selector is developed for kernel density level set estimation and its efficacy is verified in numerical studies.
\end{abstract}





\end{frontmatter}

\section{Introduction}\label{intro}

For a density function $f$ on $\mathbb{R}^d$, $d\geq 1$, its (upper) level set at a given level $c$ is defined as
\begin{align*}
\mathcal{L}_c = \{x\in\mathbb{R}^d:\;f(x) \geq c\}.
\end{align*}
With a given random sample from $f$, it is often of interest to estimate $\mathcal{L}_c$. Density level set estimation has been useful in many areas, such as clustering (Rinaldo and Wasserman, \cite{Rinaldo10}), classification (Steinwart et al., \cite{Steinwart05}), tests for multimodality (M\"{u}ller and Sawitzki, \cite{Muller91}), and topological data analysis (Fasy et al., \cite{Fasy14}).  A plug-in estimator of $\mathcal{L}_c$ using kernel density estimation is given by
\begin{align*}
\wh{\mathcal{L}}_c = \{x\in\mathbb{R}^d:\;\wh f(x) \geq c\},
\end{align*}
where $\wh f(x)$ is the kernel estimator of $f(x)$ (see (\ref{kerndens})). It is well-known that the choice of bandwidth plays a crucial role in the performance of kernel-type estimators. In this paper we derive an asymptotically optimal bandwidth for $\wh{\mathcal{L}}_c$. We take the level $c$ as a fixed value and denote $\mathcal{L}=\mathcal{L}_c$ and $\wh{\mathcal{L}}=\wh{\mathcal{L}}_c$ for simplicity. 

\begin{figure}[h]
\begin{center}
\includegraphics[scale=0.16]{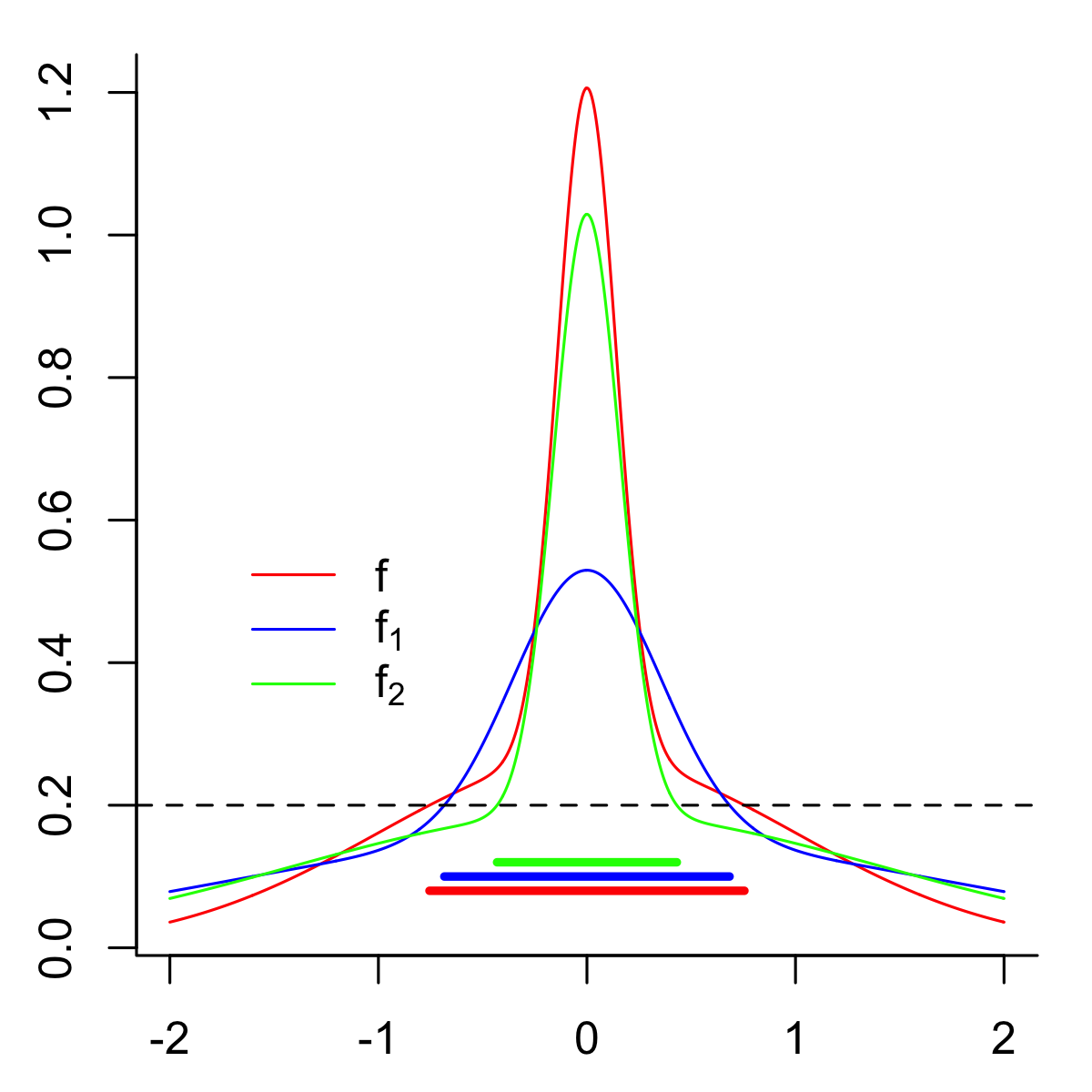}
\caption{Three density functions on $\mathbb{R}$ are shown on the graph: $f$ (red curve), $f_1$ (blue curve) and $f_2$ (green curve). The level sets with $c=0.2$ (dotted line) of the three density functions are considered. The thick horizontal lines underneath the density curves represent the level sets. Overall, the density function $f_2$ is closer than $f_1$ to $f$.  However, the level set of $f_1$ is closer to that of $f$.}
\label{fig: idea}
\end{center}
\end{figure}

The optimal bandwidth selection for $\wh f$ has been studied extensively in the literature, usually based on ISE, MISE, or MIAE (see Wand and Jones, \cite{Wand95}). These criteria measure the proximity between $\wh f$ and $f$ over $\mathbb{R}^d$. We emphasize here that the target of our estimation $\mathcal{L}$ is a set rather than a density function; this has a critical impact on the optimal bandwidth, since the quality of density estimation should be prioritized regionally rather than over the entire domain. Figure \ref{fig: idea} is an illustration, which shows that the overall closeness of density functions is not equivalent to the closeness of their level sets. Therefore, the criteria used for optimal bandwidth should be tailored specifically for nonparametric level set estimation.

A usual loss function used to measure the closeness between $\mathcal{L}$ and $\wh{\mathcal{L}}$ is based on their symmetric difference. For any two sets $A$ and $B$, let $A\Delta B$ be their symmetric difference, i.e., $A\Delta B = (A\cap B^\complement)\cup (B\cap A^\complement)$, where we use $\complement$ to denote the complement of a set. For any nonnegative integrable function $g$ on $\mathbb{R}^d$ and any Lebesgue measurable subset $A$ of $\mathbb{R}^d$, denote $\lambda_g(A) = \int_{A} g(x) dx.$ In the literature $\lambda_g (\mathcal{L}\;\Delta \; \wh{\mathcal{L}})$ has been well accepted as a measure of the proximity between $\mathcal{L}$ and $\wh{\mathcal{L}}$, due to its natural geometric interpretation. Examples of $g(x)$ include $f(x)$ and $|f(x)-c|^q$ for some $q\geq 0$. The asymptotics of $\lambda_g (\mathcal{L}\;\Delta \; \wh{\mathcal{L}})$ has been studied, e.g., in Ba\'{i}llo et al. (\cite{Baillo00}), Ba\'{i}llo (\cite{Baillo03}), Cadre (\cite{Cadre06}), Cuevas et al. (\cite{Cuevas06}), Biau et al. (\cite{Biau08}), and Mason and Polonik (\cite{Mason09}). 

We use the risk $\mathbb{E}\lambda_g (\mathcal{L}\;\Delta \; \wh{\mathcal{L}})$ to study the problem of optimal bandwidth for density level set estimation for $d\geq 1$. A critical step is to obtain the asymptotic expression for this risk, which is one of the main results in this paper and is first described below.

Let $\mathcal{M} = \{x\in\mathbb{R}^d:\;f(x) = c\}$, which is the boundary of $\mathcal{L}$ (i.e., $\mathcal{M} =\partial\mathcal{L}$) under a mild assumption (e.g., see our assumption (F1) below). Let $\text{Vol}_{d-1}$ be the natural $(d-1)$-dimensional volume measure that $\mathcal{M}$ inherits as a subset of $\mathbb{R}^d$, and $d(x,\mathcal{M})$ be the distance from any $x\in\mathbb{R}^d$ to $\mathcal{M}$. Suppose $g(x)$ is approximately $p$th power of $d(x,\mathcal{M})$ for $x$ in a small neighborhood of $\mathcal{M}$. Then under regularity conditions the following approximation holds asymptotically: 
\begin{align}\label{SymDiffL1App}
\mathbb{E}\lambda_g (\mathcal{L}\;\Delta \; \wh{\mathcal{L}}) = \mathbb{E}\int_{\mathcal{M}} |\wh f(x) - f(x)|^{p+1}w_g(x)d\text{Vol}_{d-1}(x) \{1+o(1)\},
\end{align}
where $w_g$ is a positive function on $\mathcal{M}$. Here we approximate a risk describing horizontal variations with the one constructed with vertical variations. A rigorous statement with appropriate assumptions for (\ref{SymDiffL1App}) is given in Theorem \ref{L1risk} below. Using the above expression we can interpret $\mathbb{E}\lambda_g (\mathcal{L}\;\Delta \; \wh{\mathcal{L}})$ asymptotically as a weighted $L^{p+1}$ risk for density estimation, in the form of an integration over the boundary of $\mathcal{L}$ with respect to the $(d-1)$-dimensional volume measure. When $g\equiv 1$, $\lambda_g$ is the Lebesgue measure. In this case, $p=0$ and the above approximation corresponds to the $L^1$ risk called Mean Integrated Absolute Error (MIAE), which has been used as a measure of proximity for optimal bandwidth selection for kernel density estimation (See Devroye and Gy\"{o}rfi (\cite{Devroye85}), Devroye (\cite{Devroye87}), Hall and Wand (\cite{Hall88}), Holmstr\"{o}m and Klemel\"{a} (\cite{Holmstrom92}), and Devroye and Lugosi (\cite{Devroye01})). Alternatively, Mean Integrated Squared Error (MISE) is more tractable than MIAE. This motivates us to use the choice of $g$ with $p=1$ for optimal bandwidth, specifically $g(x)=|f(x)-c|$. In this case, $\mathbb{E}\lambda_g (\mathcal{L}\;\Delta \; \wh{\mathcal{L}})$ corresponds to the excess risk (or regret) in the classification literature. In fact the excess risk has been used to find the optimal tuning parameter for nonparametric classifier, where the excess risk can be asymptotically decomposed into a squared bias term and a variance term (Hall and Kang \cite{Hall05}, Hall et al. \cite{Hall08}, Samworth \cite{Samworth12}, Cannings et al. \cite{Cannings17}). The results in this paper provide a way of understanding the excess risk as an $L^2$ risk, in a more general setting. In addition to the asymptotic approximation for the risk $\mathbb{E}\lambda_g (\mathcal{L}\;\Delta \; \wh{\mathcal{L}})$, we also show the asymptotic approximation for the error
\begin{align}\label{SymDiffL1Error}
\lambda_g (\mathcal{L}\;\Delta \; \wh{\mathcal{L}}) =  \int_{\mathcal{M}} |\wh f(x) - f(x)|^{p+1} w_g(x)d\text{Vol}_{d-1}(x) \{1+o_p(1)\},
\end{align}
under some extra assumptions on the convergence rate of the bandwidth.

Some of the important work on level set estimation includes Hartigan (\cite{Hartigan87}), Polonik (\cite{Polonik95}), Tsybakov (\cite{Tsybakov97}), Walther (\cite{Walther97}), Cadre (\cite{Cadre06}), Rigollet and Vert (\cite{Rigollet09}), among many others. Also see Mason and Polonik (\cite{Mason09}) for a comprehensive review of the literature for level set estimation. Confidence regions for level sets have recently been studied in Mammen and Polonik (\cite{Mammen13}), Sommerfeld et al. (\cite{Sommerfeld15}), Chen et al. (\cite{Chen17}), and Qiao and Polonik (\cite{QiaoPolonik19}). In Jang (\cite{Jang06}), plug-in level set estimation is applied to two-dimensional astronomical sky survey data, where the selection of bandwidth is based on the classical plug-in and cross validation approaches for density estimation. 

When the level $c$ is not explicitly given but determined by a probability value $\tau\in (0,1)$ through $c = \inf\{y\in(0,\infty): \int_{f(x)\geq y} f(x)dx\leq1-\tau\}$, we denote $c=c(\tau)$ and $\mathcal{L}_{c(\tau)}$ is called the $100(1-\tau)\%$ highest density region (HDR) of $f$ (see Hyndman, \cite{Hyndman96}). The corresponding plug-in estimator is $\wh{\mathcal{L}}_{\wh c(\tau)}$, where  $\wh c(\tau) =\inf\{y\in(0,\infty): \int_{\wh f(x)\geq y} \wh f(x)dx \leq 1- \tau\}$. In the case of $f$ being a univariate density (i.e., $d=1$), the bandwidth selection problem for estimating HDR was studied in Samworth and Wand (\cite{Samworth10}). They chose $\mathbb{E} \lambda_g (\mathcal{L}_{c(\tau)}\;\Delta \; \wh{\mathcal{L}}_{\wh c(\tau)})$ with $g=f$ as the risk function to minimize for bandwidth selection. The extension of their approach to the multivariate case is far from trivial and has been recently studied in Doss and Weng (\cite{Doss18}).

The work in Doss and Weng (\cite{Doss18}) also considers the bandwidth selection problem of the estimation of density level sets. The comparison of their work with an earlier arXiv version of the present paper (see Qiao, \cite{Qiao18}) has been discussed in Doss and Weng (\cite{Doss18}). In particular, the risk criterion they use for bandwidth selection for level set estimation is $\mathbb{E}\lambda_g (\mathcal{L}\;\Delta \; \wh{\mathcal{L}})$ with $g=f$, which is an $L^1$ type of risk as a special case of (\ref{SymDiffL1App}). See Remark~\ref{cororemark} for Corollary~\ref{L1asympexpress} as well as the Discussion section for more detailed comparisons.

The rest of the paper is organized as follows. We first introduce some notation and geometric concepts in Section 2. In Section 3, after discussing the assumptions that we will use, we derive some asymptotic results for the $L^p$ type of risks introduced above and an optimal bandwidth for density level set estimation. Specifically, Theorems \ref{L1risk} and \ref{riskexpect} formulate the ideas given in (\ref{SymDiffL1Error}) and (\ref{SymDiffL1App}), respectively. Corollary \ref{L1asympexpress} gives an exact asymptotic expression for $\mathbb{E}\lambda_g (\mathcal{L}\;\Delta \; \wh{\mathcal{L}})$ when $p=0$. The excess risk as an asymptotic $L^2$ type of risk, is used to find the optimal bandwidth with the result given in Theorem \ref{optimAMISE}. Simulation results are presented in Section 4, where we show the efficacy of our bandwidth selector in finite samples. We leave all the proofs to Section 5, while some miscellaneous results are put in the appendix.

\section{Notation and some geometric concepts}\label{notationsec}

Let $X_1,\cdots, X_n$ be i.i.d. from the $d$-dimensional density function $f$. Denote the bandwidth vector ${\bf h}=(h_1,h_2,\cdots,h_d)^T$ and ${\bf h}^{-1}=(h_1^{-1},h_2^{-1},\cdots,h_d^{-1})^T$. We consider a kernel density estimator 
\begin{align}\label{kerndens}
\wh f(x) = \frac{1}{n\Pi_{j=1}^dh_j } \sum_{i=1}^n K\left({\bf h}^{-1} \odot (x-X_i)\right),
\end{align}
where $K$ is a kernel function on $\mathbb{R}^d$ and $\odot$ is used to denote the Hadamard or element-wise product between two vectors of the same size. Here we assign a bandwidth value for each of the variables in the density estimation. This corresponds to a diagonal bandwidth matrix, which is a compromise between flexibility (by using a full bandwidth matrix) and simplicity (by using only a scalar bandwidth). See Wand and Jones (\cite{Wand94}) for more discussion on the impact of the form of bandwidth matrix on multivariate density estimation. We use a product kernel for $K$, i.e., we can write
\begin{align*}
K\left({\bf h}^{-1} \odot (x-X_i)\right) = \prod_{j=1}^d \wt K\left( \frac{x_j-X_{ij}}{h_j}\right),
\end{align*}
where $\wt K$ is a univariate kernel function, $X_i=(X_{i1},\cdots,X_{id})^T$, for $i=1,\cdots,n$, and $x=(x_1,\cdots,x_d)^T$. The order of a kernel is determined by its first nonzero moment. We call $K$ a $\nu$th ($\nu\geq2$) order kernel if $\int_{\mathbb{R}} |u^\nu\wt K(u)| du<\infty$ and
\begin{align*}
    \int_{\mathbb{R}} u^l \wt K(u) du=
    \begin{cases}
      1, & \text{if } l=0, \\
      0, & \text{if } l=1,\cdots,\nu-1,\\
      \kappa_\nu \neq 0, & \text{if } l=\nu.
    \end{cases}
%
%
\end{align*}
It is obvious that $\nu$ is always even if $\wt K$ is symmetric. When $d=1$, we also denote $h=\mathbf{h}$ and write $\wh f(x)=\frac{1}{nh}\sum_{i=1}^n K(h^{-1}(x-X_i))$.

{\em Notation}:  Let $\mathscr{H}_{d-1}$ be the $(d-1)$-dimensional normalized Hausdorff measure on $\mathbb{R}^d$ (cf. Evans and Gariepy, \cite{Evans92}). It agrees with the $(d-1)$-dimensional volume measure $\text{Vol}_{d-1}$ on nice sets. For $d=1$, $\mathscr{H}_0$ is the cardinality of a set, such that for $A=\{a_1,\cdots, a_m\}\subset\mathbb{R}$ and a function $g:\mathbb{R}\mapsto\mathbb{R}$, $\int_A g(x) d\mathscr{H}_0(x) = \sum_{i=1}^m g(a_i)$. For simplicity, we usually omit the subscript $d-1$ of $\mathscr{H}_{d-1}$ and write $\mathscr{H} = \mathscr{H}_{d-1}$  if the value of $d$ is clear in the context. Let $\lambda$ be the $d$-dimensional Lebesgue measure ($d\geq 1$). Recall that $\lambda$ and $\mathscr{H}_{d}$ are equal on $\mathbb{R}^d$. Let $\mathbf{I}$ be the indicator function. 

For a positive integer $p$, let $\mathbb{Z}_+^p$ be the set of all the $p$-dimensional vectors of positive integers. For any $\mathbf{i}=(i_1,\cdots,i_p)^T\in\mathbb{Z}_+^p$, let $\max(\mathbf{i}) = \max(i_1, \cdots, i_p)$. Denote $\mathbb{Z}_+^{p,d}=\{\mathbf{i}\in\mathbb{Z}_+^p,\max(\mathbf{i})\leq d\}$. For any function $g: \mathbb{R}^d\mapsto \mathbb{R}$ with $p$th (partial) derivatives ($p\geq 1$), and for $\mathbf{i}\in\mathbb{Z}_+^{p,d}$,  denote $g_{(\mathbf{i})}(x) = g_{(i_1,\cdots,i_p)}(x) = \frac{\partial^p}{\partial x_{i_1}\cdots \partial x_{i_p}} g(x)$ with the convention $g_{\mathbf{(i)}}(x) = g(x)$ for $\mathbf{i}\in\mathbb{Z}_+^0$. For example,  $g_{(k,l)}(x) = \frac{\partial^2}{\partial x_k \partial x_l} g(x)$ for $1 \leq k,l \leq d$. If $i_1=\cdots=i_p =i$ for $1\leq i\leq d$, we denote $g_{(i*p)}(x)=g_{(i_1,\cdots,i_p)}(x)$. For $x\in\mathbb{R}^d$ and $\mathbf{i}\in\mathbb{Z}_+^{p,d}$, denote $x^{(\mathbf{i})}=x_{i_1}\times\cdots \times x_{i_p}.$ 
For $d\geq 2$, we denote the gradient and Hessian matrix of $g$ by $\nabla g$ and $\nabla^2 g$, respectively. With slight abuse of notation, we also use $\nabla g$ to denote the first derivative $g^\prime$ when $d=1$. For any Borel set $A\subset\mathbb{R}$, let $g^{-1}(A) = \{x\in\mathbb{R}^d:\; g(x)\in A\}$. Let $\|g\|_q = (\int_{\mathbb{R}^d} |g(x)|^qdx)^{1/q}$ for $q>0$ and $\|g\|_\infty = \sup_{x\in\mathbb{R}^d} |g(x)|$. For sequences $a_n, b_n\in\mathbb{R}$, we denote $a_n \asymp b_n$ if $0<\lim\inf_{n\rightarrow\infty}|a_{n}/b_n|\leq \lim\sup_{n\rightarrow\infty}|a_{n}/b_n| < \infty$. For a sequence $\mathbf{a}_n=(a_{1,n},\cdots,a_{d,n})^T\in\mathbb{R}^d$, $d\geq2$, denote $\mathbf{a}_n\asymp b_n$ if $a_{i,n}\asymp b_n$, $i=1,\cdots,d$.

Next we introduce some geometric concepts. For any set $A\subset\mathbb{R}^d$ and $\epsilon>0$, we denote the $\epsilon$-enlargement of $A$ by $A\oplus\epsilon = \bigcup_{x\in A} \mathcal{B}_x(\epsilon)$, where $\mathcal{B}_x(\epsilon) = \{y\in\mathbb{R}^d:\; \|x-y\|\leq \epsilon\}.$ For any two sets $A, B\subset\mathbb{R}^d$, let $d_H(A, B)$ be the Hausdorff distance between $A$ and $B$, i.e., $$d_H(A, B) = \max\left\{\sup_{x\in B} d(x,A),\; \sup_{x\in A}d(x,B)\right\},$$ where $d(x,A)=\inf_{y\in A}\|x-y\|$. Let $\pi_A(x)$ be the set of the closest points in $A$ to $x$, i.e. $\pi_A(x) =\{y\in A:\; \|x-y\|=d(x,A)\}$, which is called the normal projection of $x$ onto $A$. 

\begin{figure}[h]
\begin{center}
\includegraphics[scale=0.26]{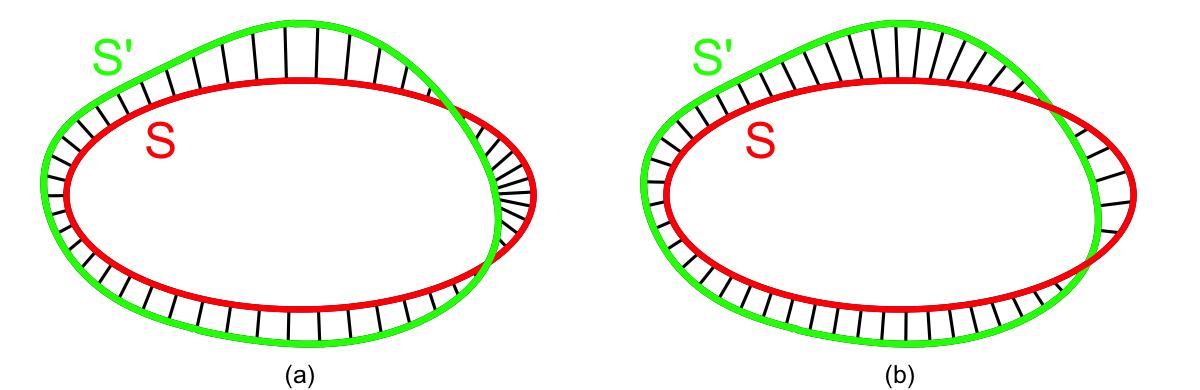}
\caption{Two normal compatible curves $S$ and $S^\prime$. (a) represents the normal projection from $S^\prime$ to $S$. (b) represents the normal projection from $S$ to $S^\prime$.}
\label{fig: normalcomp}
\end{center}
\end{figure}

Let $A$ and $B$ be two $(d-1)$-dimensional smooth submanifolds embedded in $\mathbb{R}^d$ ($d\geq 2$). Then the normal projections $\pi_A: B \mapsto A$ and $\pi_B: A\mapsto B$ define two maps between $A$ and $B$. The two manifolds $A$ and $B$ are called {\em normal compatible} if the projections $\pi_A$ and $\pi_B$ are homeomorphisms. See Chazal et al. (\cite{Chazal07}) and Figure \ref{fig: normalcomp} for a graphical illustration.

We will also use the concept of {\em reach} of a manifold. For a $p$-dimensional manifold $\mathcal{S}$ embedded in $\mathbb{R}^d$ ($p<d$), the reach of $\mathcal{S}$, denoted by $\rho(\mathcal{S})$, is the largest $\delta$ such that the normal projection from every point in $\mathcal{S}\oplus\delta$ onto $\mathcal{S}$ is unique. See Federer (\cite{Federer59}). A positive reach corresponds to the notion of bounded curvature of a manifold. See Niyogi et al. (\cite{Niyogi08}) and Genovese et al. (\cite{Genovese12}).

\section{Main results}

\subsection{Assumptions and their discussion}

We introduce the assumptions that will be used in this paper. Let $\mathbb{R}_+^d$ be the set of vectors in $\mathbb{R}^d$ with positive coordinates. With the requirement of $C^2$ smoothness of the kernel function $K$, define the class of functions
\begin{align*}
\mathscr{K} =\left\{ K_{(\mathbf{i})}(\mathbf{h}^{-1}\odot(x-\cdot)): \mathbf{h}\in\mathbb{R}_+^d, x\in\mathbb{R}^d,\mathbf{i}\in\mathbb{Z}_+^p, p\in\{0,1,2\}, \max(\mathbf{i})\leq d\right\}.
\end{align*}
Let $\mathscr{B}$ be the Borel $\sigma$-algebra on $\mathbb{R}^d$. For any probability measure $Q$ on $(\mathbb{R}^d, \mathscr{B})$ and $\epsilon>0$, let $N(\mathscr{K}, L_2(Q), \epsilon)$ be the $\epsilon$-covering number for $\mathscr{K}$ using the $L^2$ norm with respect to $Q$, i.e., the minimal number of balls $\{g: \int_{\mathbb{R}^d} |g-\wt g|^2 dQ \leq \epsilon\}$ needed to cover $\mathscr{K}$. Denote the envelope function $F(x)=\sup_{g\in\mathscr{K}}|g(x)|$. For $\delta>0$, denote $\mathcal{I}(\delta) = f^{-1}([c-\delta/2,c+\delta/2])$. For any $x\in\mathbb{R}^d$ such that $\nabla f(x)\neq 0$ and $s\in\mathbb{R}$, denote $$\zeta_x(s) = x+\frac{\nabla f(x)}{\|\nabla f(x)\|}s.$$ Note that $\zeta_x(s) =x+ \text{sign}(f^\prime(x))\times s$ when $d=1$. For $1\leq p\leq\infty$, denote the $L^p$ space of Lebesgue measurable functions on $\mathbb{R}^d$ by $\mathscr{L}^p$.\\ [-8pt]

\hspace{-12pt}{\bf Assumptions}: \\[-10pt]

{\bf (F1)} $f$ is a $\nu$ times continuously differentiable pdf for some $\nu\geq2$. The density function and all of its first to $\nu$-th derivatives are bounded on $\mathbb{R}^d$. We also assume $C_{\ell}<c<C_u$, where $C_{\ell}=\inf_{x\in\mathbb{R}^d} f(x)$ and $C_u=\sup_{x\in\mathbb{R}^d} f(x)$, and that there exist $\delta_0>0$ and $\epsilon_0>0$ such that  $\|\nabla f(x)\| > \epsilon_0$ for all $x\in \mathcal{I}(2\delta_0)$.

{\bf (G1)} $g$ is a non-negative continuous function on $\mathbb{R}^d$ and there exist $p\geq0$ and a bounded positive function $g^{(p)}(x)$ on $\mathcal{M}$ such that as $s\rightarrow 0$,
\begin{align*}
\sup_{x\in\mathcal{M}} \left|\frac{g(\zeta_x(s))}{|s|^p} - g^{(p)}(x)\right| = o(1).
\end{align*}

{\bf (K1)} $K$ is a symmetric product kernel function of $\nu$th order for some $\nu\geq 2$. Also $K\in\mathscr{L}^1\cap\mathscr{L}^\infty$. 

{\bf (K2)} $K$ is two times continuously differentiable. We require that $\|F\|_\infty<\infty$ and for some $C_0>0$ and $\eta>0$, $$\sup_Q N(\mathscr{K}, L_2(Q), \epsilon \|F\|_\infty)\leq C_0\epsilon^{-\eta},$$ for $0<\epsilon<1$, where the supremum is taken over all the probability measures $Q$ on $(\mathbb{R}^d, \mathscr{B})$.\\[-8pt]

%


\begin{remark}\label{assumpremark}$\;$\\[-14pt]
{\normalfont
\begin{itemize}
\item[a)]  In assumption (F1), the global smoothness requirement for $f$ can be weakened to only hold regionally on $\mathcal{I}(2\delta_0)$, if we choose to use a kernel function $K$ with bounded support. Conditions similar to $\|\nabla f(x)\| > \epsilon_0$ for $x\in\mathcal{I}(2\delta_0)$ in assumption (F1) have appeared in Cadre (\cite{Cadre06}), Cuevas et al. (\cite{Cuevas06}), Mammen and Polonik (\cite{Mammen13}), among others. It excludes the possibility of ``flat parts'' around the level set. In particular, it implies that $\mathcal{M} = \partial\mathcal{L}$, which is a compact $(d-1)$-dimensional $C^1$ submanifold in $\mathbb{R}^d$ (see Theorem 2 in Walther (\cite{Walther97})). In the case $d=1$, $\mathcal{M}$ is a collection of separated points, i.e., there exist $x_1,\cdots,x_N$ for some positive integer $N$ such that $\mathcal{M} = \{x_i:\; i=1,2,\cdots,N\}$.

\item[b)] An assumption similar to (G1) has appeared in Mason and Polonik (\cite{Mason09}). Below we give the specific forms of $g^{(p)}$ for some usual functions $g$.
\begin{itemize}
\item[(i)] If $g$ is a continuous function with positive values on $\mathcal{M}$, then $p=0$ and $g^{(p)}(x) =g(x)$, $x\in\mathcal{M}$. Examples include $g(x)\equiv1$ and $g(x) = f(x)$. 
%
\item[(ii)] If $g(x) = f(x)^r|f(x)-c|^q$ for some $q>0$ and any $r\geq 0$,  then $p =q$ and $g^{(q)}(x) =c^r\|\nabla f(x)\|^q$,  $x\in\mathcal{M}$. This is because for $x\in\mathcal{M}$, as $s\rightarrow 0$,
\begin{align*}
&g(\zeta_x(s)) \\
=& [f(x+s\times\nabla f(x)/\|\nabla f(x)\|)]^r \times |f(x+s\times\nabla f(x)/\|\nabla f(x)\|)-c|^q \\
= &[c+o(s)]^r \times \Large|s\|\nabla f(x)\| + o(s)\Large|^q \\
= & c^r |s|^q \|\nabla f(x)\|^q +o(|s|^q).
\end{align*}
\end{itemize}

\item[c)] For assumption (K1), it is known that using higher order kernels (i.e. $\nu>2$), together with higher order smoothness assumptions can reduce the bias in kernel density estimation. But higher order kernels are avoided sometimes because it is possible that density estimators have negative values (see, e.g., Silverman, \cite{Silverman86}, page 69). However, this should be of less concern for level set estimation, because the negative values of the density estimate are not (directly) involved in our level set estimator for $c> 0$. 

\item[d)] Assumption (K2) is imposed to uniformly control the stochastic variation of the kernel density estimator and their derivatives around the expectations. Similar assumptions have appeared in Gin\'{e} and Guillou (\cite{Gine02}), and Einmahl and Mason (\cite{Einmahl05}). Also see Chen et al. (\cite{Chen17}). For sufficient conditions for (K2) to hold, see e.g., Nolan and Pollard (\cite{Nolan87}). In particular, the Gaussian kernel and many usual kernels with bounded support satisfy assumption (K2).
\end{itemize}
}
\end{remark}

\subsection{Asymptotic expressions of $\mathbb{E}\lambda_g (\mathcal{L}\;\Delta \; \wh{\mathcal{L}})$ and $\lambda_g (\mathcal{L}\;\Delta \; \wh{\mathcal{L}})$}

For $x\in\mathcal{M}$, let 
\begin{align*}
t_n(x) = \argmin_{t}\left\{|t|:\; \zeta_x( t ) \in \wh{\mathcal{M}} \right\}.
\end{align*}
For $d\geq 2$, once we establish the normal compatibility between $\wh{\mathcal{M}}$ and $\mathcal{M}$, the inverse mapping of the normal projection $\pi_{\mathcal{M}}$ will be well defined and is denoted by $P_n$. Namely, for any $x\in\mathcal{M}$, we have $P_n(x)\in\wh{\mathcal{M}}$ and $P_n(x)-x$ is orthogonal to the tangent space of $\mathcal{M}$ at $x$. We also write $\wh{\cal{M}} = P_n(\cal{M})$. Since $ \nabla f(x)$ is a normal vector of $\mathcal{M}$ at $x$, we can write 
\begin{align}\label{map}
P_n(x) = \zeta_x( t_n(x)),
\end{align}
for some unique  $t_n(x)\in\mathbb{R}$. For $d=1$, $P_n$ is set to be equivalent to $\pi_{\wh{\mathcal{M}}}$, i.e. it maps points in $\mathcal{M}$ to their closest points in $\wh{\mathcal{M}}$. 


Let $\{A_i: \;i=1,2,\cdots,N_d\}$ be a partition of $\mathcal{M}$ and $a_i$ be a point on $A_i$.  Since pointwisely $t_n$ is small when $n$ is large, the following approximation is heuristic when the partition is fine enough:
\begin{align*}
\lambda (\mathcal{L}\;\Delta \; \wh{\mathcal{L}}) \approx \sum_{i=1}^{N_d} |t_n(a_i)|\mathscr{H}(A_i) \approx \int_{\mathcal{M}} |t_n(x)|d\mathscr{H}(x). 
\end{align*}
The more precise form of the above idea is given in the following theorem, where we need the assumption on $n$ and $\mathbf{h}$ as below.
\begin{itemize}
\item[]{\bf (H1)} The bandwidth (vector) $\mathbf{h}\in\mathbb{R}_+^d$ is dependent on $n$ such that 
\begin{align*}
(\log n)^{-1}nh_1\cdots h_d\|\mathbf{h}\|^{4}\rightarrow\infty \text{ and } \log(1/\|\mathbf{h}\|)/(\log\log{n})\rightarrow \infty,
\end{align*}
as $n\rightarrow\infty$. When $d\geq2$, we assume $h_i \asymp h_j$, for $1\leq i,j \leq d$.
\end{itemize}
\begin{theorem}\label{L1risk}
Under assumptions (K1), (K2), (F1), (G1) and (H1), we have
\begin{align}\label{L1express}
\lambda_g (\mathcal{L}\;\Delta \; \wh{\mathcal{L}}) =  \frac{1}{p+1} \int_{\mathcal{M}}\frac{g^{(p)}(x)}{\|\nabla f(x)\|^{p+1}} |\wh f(x) - f(x)|^{p+1}d\mathscr{H}(x) \{1+o_p(1)\}.
\end{align}
\end{theorem}

\begin{remark} $\;$\\[-14pt]
{\normalfont
\begin{itemize}
\item[a)]  This result is related to but different from Theorem 2.1 in Cadre (\cite{Cadre06}), where assumptions are imposed to ensure that $\sqrt{nh^d} \lambda_g (\mathcal{L}\;\Delta \; \wh{\mathcal{L}}) \rightarrow \mu_g$ in probability for some $\mu_g>0$ as $n\rightarrow\infty$. In particular the bandwidth is assumed to be small enough that the bias in the kernel density estimation can be ignored. In contrast, our focus is on revealing the asymptotic expression of $\lambda_g (\mathcal{L}\;\Delta \; \wh{\mathcal{L}})$ for the purpose of finding the optimal bandwidth, for which both the variance and bias in the kernel density estimation are involved.
\item[b)]  The assumption $(\log n)^{-1}nh_1\cdots h_d\|\mathbf{h}\|^{4}\rightarrow\infty$ in this theorem is used to guarantee the normal compatibility between $\mathcal{M}$ and $\wh{\mathcal{M}}$ for $d\geq2$, and can in fact be relaxed and replaced with $(\log n)^{-1}nh^3\rightarrow\infty$ for $d=1$, which is required for the uniform consistency of the kernel estimation for the first derivative of the density. As indicated in Section~\ref{notationsec}, $\mathscr{H}_0$ is the cardinality measure. For $d=1$, with $\mathcal{M} = \{x_i:\; i=1,2,\cdots,N\}$ (see the discussion after the assumptions), the result (\ref{L1express}) becomes 
\begin{align*}
 \lambda_g (\mathcal{L}\;\Delta \; \wh{\mathcal{L}}) = \frac{1}{p+1} \sum_{i=1}^N \frac{|\wh f(x_i) - f(x_i)|^{p+1}}{|f^\prime(x_i)|^{p+1}}g^{(p)}(x_i)\{1+o_p(1)\}.
\end{align*}
\end{itemize}
}
\end{remark}

The required assumption of $(\log n)^{-1}nh_1\cdots h_d\|\mathbf{h}\|^{4}\rightarrow\infty$ is critical in the above theorem. However, if we only consider the expectation of $\lambda_g (\mathcal{L}\;\Delta \; \wh{\mathcal{L}})$, it is in fact not needed. We modify (H1) into the following weaker assumption.
\begin{itemize}
\item[] {\bf(H2)} The bandwidth (vector) $\mathbf{h}\in\mathbb{R}_+^d$ is dependent on $n$ such that 
\begin{align*}
(\log n)^{-1}nh_1\cdots h_d\rightarrow\infty \text{ and } \log(1/\|\mathbf{h}\|)/(\log\log{n})\rightarrow \infty,
\end{align*}
as $n\rightarrow\infty$. When $p\geq4$ where $p$ appears in assumption (G1), we further assume that $(\log n)^{-(p-2)}nh_1\cdots h_d\rightarrow\infty$. When $d\geq2$, we assume $h_i \asymp h_j$, for $1\leq i,j \leq d$.\\[-10pt]
\end{itemize}
 Let $s_n>0$ be such that
 \begin{align}
 s_n^2 &= \frac{1}{nh_1\cdots h_d}\| K\|_2^2 c,\label{snsq}\\
\text{and } \beta_{\mathbf{h}}(x) &= \frac{1}{\nu!}\kappa_\nu\sum_{k=1}^d h_k^\nu f_{(k * \nu)}(x).\label{betahx}
 \end{align}
Notice that $\beta_{\mathbf{h}}(x)=O(\|\mathbf{h}\|^\nu)$ if the boundedness of the $\nu$-th derivatives of $f$ is assumed (see F1). It is known (see, e.g., Wand and Jones, \cite{Wand95}; also see (\ref{fact2}) and (\ref{fact1}) in the proof) that under regularity conditions the bias for kernel density estimator at $x\in\mathcal{M}$ is
\begin{align}\label{MISEbias}
\mathbb{E} \wh f(x) -f(x) = \beta_{\mathbf{h}}(x) + o(\|\mathbf{h}\|^\nu),
\end{align}
and variance is
\begin{align}\label{MISEvar}
\text{Var}(\wh f(x)) = s_n^2(1+o(1)).
\end{align}
We have the following theorem. 

\begin{theorem}\label{riskexpect}
Under assumptions (K1), (K2), (F1), (G1), (H2), we have
\begin{align}
&\mathbb{E}\lambda_g (\mathcal{L}\;\Delta \; \wh{\mathcal{L}}) \nonumber\\
=  & \frac{1}{1+p} \mathbb{E} \int_{\mathcal{M} } \frac{g^{(p)}(x) }{\|\nabla f(x)\|^{p+1}} \left| s_n\,Z + \beta_{\mathbf{h}}(x) \right|^{p+1} d\mathscr{H}(x) + o(s_n^{p+1} + \|\mathbf{h}\|^{\nu(p+1)}), \label{L1expexpress}\\
=  & \frac{1}{1+p} \mathbb{E}\int_{\mathcal{M}}\frac{g^{(p)}(x)}{\|\nabla f(x)\|^{p+1}} |\wh f(x) - f(x)|^{p+1}d\mathscr{H}(x) + o(s_n^{p+1} + \|\mathbf{h}\|^{\nu(p+1)}). \label{L1express2}
\end{align}
where $Z$ is a standard normal random variable. 
\end{theorem}

\begin{remark}{\normalfont
Using the symmetry of $Z$'s distribution, we have
\begin{align*}
\mathbb{E} \left| s_n\,Z + \beta_{\mathbf{h}}(x) \right|^{p+1} = \mathbb{E} \left| s_n\,Z + |\beta_{\mathbf{h}}(x)| \right|^{p+1} \geq \max[s_n^{p+1} \mathbb{E}(|Z|^{p+1}), |\beta_{\mathbf{h}}(x)|^{p+1}].
\end{align*}
Also see (\ref{normlowerbound}) in the proof for a lower bound. So the first terms on the right-hand sides of (\ref{L1expexpress}) and (\ref{L1express2}) are indeed leading terms, if $|\beta_{\mathbf{h}}(x)|/\|\mathbf{h}\|^\nu$ is not zero for all $x\in\mathcal{M}$.
}
\end{remark}

Notice that $g\equiv g^{(p)}$ when $p=0$ in assumption (G1). By observing the fact for any $a\in\mathbb{R}$,
\begin{align}
\mathbb{E}|Z-a| = |a| \mathbb{P}(|Z|\leq |a|) + \sqrt{\frac{2}{\pi}} e^{-a^2/2} = \gamma(|a|),
\end{align}
where 
\begin{align*}
\gamma(u) &= \sqrt{\frac{2}{\pi}} \left( u\int_0^u e^{-t^2/2} dt + e^{-u^2/2}\right), \; u\geq 0,
\end{align*}
we have the following corollary which gives an exact asymptotic expression of $\mathbb{E}\lambda_g (\mathcal{L}\;\Delta \; \wh{\mathcal{L}})$ when $p=0$, the example including $g\equiv1$ and $g=f$. The result is comparable to Theorem 1 in Devroye and Gy\"{o}rfi (\cite{Devroye85}, page 78), where they considered the MIAE as the risk for kernel density estimation. 
\begin{corollary}\label{L1asympexpress}
Suppose $p=0$ in assumption (G1). Under assumptions (K1), (K2), (F1), (G1), (H2), we have
\begin{align}\label{L1exact}
 \mathbb{E} \lambda_g (\mathcal{L}\;\Delta \; \wh{\mathcal{L}}) = \int_{\mathcal{M}} \frac{s_n\gamma(|\beta_{\mathbf{h}}(x)|/s_n)}{\|\nabla f(x)\|} g(x)d\mathscr{H}(x) + o \left(\|{\bf h}\|^\nu + \frac{1}{\sqrt{nh_1\cdots h_d}}\right).
\end{align}

\end{corollary} 
\begin{remark} \label{cororemark}$\;$\\[-10pt]
{\normalfont
\begin{itemize}
\item[a)]  One can obtain asymptotic lower and upper bounds for the risk in (\ref{L1exact}), following similar arguments as in the proof of Theorem 2 in Devroye and Gy\"{o}rfi (\cite{Devroye85}, page 79), or in Holmstr\"{o}m and Klemel\"{a} (\cite{Holmstrom92}, page 257). For example, for the upper bound, since $\gamma(u)\leq u + \sqrt{2/\pi}$, $u\geq 0$, we have
\begin{align}\label{upperbound}
 &\mathbb{E} \lambda_g (\mathcal{L}\;\Delta \; \wh{\mathcal{L}}) \nonumber\\
 \leq &\left[  \int_{\mathcal{M}} \frac{\left|\beta_{\mathbf{h}}(x)\right| g(x)}{\|\nabla f(x)\|} d\mathscr{H}(x) + \sqrt{\frac{2}{\pi}} s_n \int_{\mathcal{M}} \frac{g(x)}{\|\nabla f(x)\|} d\mathscr{H}(x)\right] + o(\|{\bf h}\|^\nu +s_n).
\end{align}
The minimization of the upper bound leads to an approximation to the asymptotically optimal bandwidth, as the closed form of the minimizer for the leading term in (\ref{L1exact}) is difficult to obtain. See Devroye and Gy\"{o}rfi (\cite{Devroye85}, page 107) for a similar suggestion. In the case $h_1=\cdots = h_d=h$, the leading term in the above upper bound can be analytically minimized with respect to $h$. Using a numerical method following the ideas in Hall and Wand (\cite{Hall88}), where minimizing the MIAE of kernel density estimation is considered, it is also possible to find an asymptotic optimal bandwidth selector tailored for the level set estimation by minimizing $\mathbb{E}\int_{\mathcal{M}} \frac{|\wh f(x) - f(x)|}{\|\nabla f(x)\|} d\mathscr{H}(x).$
\item[b)] If we specifically choose $g=f$ and $\nu=2$, then the result in this corollary is similar to Theorem 2.1 in Doss and Weng (\cite{Doss18}), where they consider the selection of a bandwidth matrix for level set estimation. In fact Theorem 2.1 in Doss and Weng (\cite{Doss18}) can be understood as a special case of this corollary, if their bandwidth matrix is restricted to be diagonal. In this corollary we approximate $\mathbb{E} \lambda_f (\mathcal{L}\;\Delta \; \wh{\mathcal{L}})$ as an $L^1$ type of risk, which is a special case of a more general result in Theorem~\ref{riskexpect}.
\end{itemize}
}
\end{remark}

In addition to the case $p=0$ covered in Corollary~\ref{L1asympexpress}, another interesting scenario is $p=1$ in  assumption (G1), which holds when $g(x)=g_r(x):= f(x)^r|f(x)-c|$ for some $r\geq 0$. Note that the choice of $r$ only impacts up to a constant in the asymptotic form of $\mathbb{E}\lambda_{g}(\mathcal{L}\;\Delta \; \wh{\mathcal{L}})$ in Theorem~\ref{riskexpect} when $g=g_r$. This is because $g_r^{(1)}(x) = c^r \|\nabla f(x)\|$ (see the calculation in Remark~\ref{assumpremark} b)(ii)). 

We call the quantity $\mathbb{E}\lambda_{g_r}(\mathcal{L}\;\Delta \; \wh{\mathcal{L}})$ the ``excess risk'', for $r\geq 0$. This is closely related to the concept of excess risk frequently used in the classification literature (see, e.g. Samworth, \cite{Samworth12}). Suppose we have a random pair $(X,Y)\in \mathbb{R}^d \times \{0,1\}$, where $Y$ is the class label of $X$. Then the Bayes optimal classifier is $\psi(x):=\mathbf{I}(\eta(x)\geq \frac{1}{2})$, where $\eta(x)=\mathbb{E}(Y|X=x)$. The misclassification risk is $R(\psi) : = \mathbb{P}(\psi(X)\neq Y)$. Given an i.i.d. sample $\mathcal{X}$ with the same distribution as $(X,Y)$, suppose one can find an estimator $\wh\eta(x)$ for $\eta(x)$ and build an empirical classifier $\wh\psi(x):=\mathbf{I}(\wh\eta(x)\geq \frac{1}{2})$. Then the misclassification risk for $\wh\psi$ is $\mathbb{E}R(\wh\psi) = \mathbb{P}(\wh\psi(X)\neq Y)$. The difference $\mathbb{E}R(\wh\psi) - R(\psi)$ is called the excess risk in this binary classification problem and it is well-known that one can write 
\begin{align}\label{regret}
\mathbb{E}R(\wh\psi) - R(\psi) = \frac{1}{2} \mathbb{E}\int_{\mathcal{S}_{\wh\eta}\Delta \mathcal{S}_{\eta}} \left|\eta(x)-\frac{1}{2}\right|dx,
\end{align}
where $\mathcal{S}_{\eta} = \{x\in\mathbb{R}^d: \;\eta(x)\geq \frac{1}{2}\}$ and $\mathcal{S}_{\wh\eta}=\{x\in\mathbb{R}^d: \; \wh\eta(x)\geq \frac{1}{2}\}$ are level sets of $\eta$ and $\wh\eta$ at the level $\frac{1}{2}$, respectively. Note that the above expression has the same form as $\mathbb{E}\lambda_{g_0}(\mathcal{L}\;\Delta \; \wh{\mathcal{L}})$.

\subsection{Optimal bandwidth using excess risk}\label{L2TyprRisk}

We use the excess risk defined above to find the asymptotic optimal bandwidth for density level set estimation, motivated by the connection to the classification literature. In fact, the excess risk is also studied in the literature of density level set estimation (see, e.g. Rinaldo and Wasserman, \cite{Rinaldo10}). For any measurable set $\mathcal{A}\subset \mathbb{R}^d$, the {\em excess mass functional} is defined as $\mathcal{E}(\mathcal{A}) =\mathbb{P}(\mathcal{A})  - c\lambda(\mathcal{A})$. It is known that $\mathcal{E}(\mathcal{A})$ is maximized when $\mathcal{A}=\mathcal{L}$. Then it is easy to show that $\mathcal{E}(\mathcal{L}) - \mathbb{E}[\mathcal{E}(\wh{\mathcal{L}})] = \mathbb{E}\lambda_{g_0}(\mathcal{L}\;\Delta \; \wh{\mathcal{L}})$. In other words, minimizing the risk $\mathbb{E}\lambda_{g_0}(\mathcal{L}\;\Delta \; \wh{\mathcal{L}})$ is equivalent to maximizing $\mathbb{E}[\mathcal{E}(\wh{\mathcal{L}})]$. The excess risk $\mathbb{E}\lambda_{g_0}(\mathcal{L}\;\Delta \; \wh{\mathcal{L}})$ is ``cost-sensitive'' (Scott and Davenport, \cite{Scott06}) in the sense that the weight function $g_0(x)=|f(x)-c|$ penalizes more heavily at a point $x\in \mathcal{L}\;\Delta \; \wh{\mathcal{L}}$, if its density value deviates more from the level $c$.

One can also understand the excess risk for density level set estimation from a binary classification perspective. Given a random vector $X\sim f$, which is independent of $X_1,\cdots,X_n$, suppose that we would like to find a set $\mathcal{A}$, such that we claim $f(X)\geq c$ when $X\in \mathcal{A}$ and $f(X)<c$ when $X\in {\mathcal{A}}^\complement$, where ${\mathcal{A}}^\complement$ is the complement of $\mathcal{A}$. Define the loss function 
\begin{align}
e_{\mathcal{A}}(x) = [c-f(x)] [\mathbf{I}(x\in \mathcal{A}) - \mathbf{I}(x\in {\mathcal{A}}^\complement)],\; x\in\mathbb{R}^d.
\end{align}
Notice that this loss function is related to the excess mass functional through $\int_{\mathbb{R}^d}e_{\mathcal{A}}(x) dx = \mathcal{E}({\mathcal{A}}^\complement) - \mathcal{E}(\mathcal{A}).$ Also it is clear that $\mathcal{A}=\mathcal{L}$ minimizes the risk function $\mathcal{R}(\mathcal{A}):=\mathbb{E}[e_{\mathcal{A}}(X) ]$. Note that $$\mathbb{E}\mathcal{R}(\wh{\mathcal{L}}) - \mathcal{R}(\mathcal{L}) = 2 \mathbb{E}\lambda_{g_1}(\mathcal{L}\;\Delta \; \wh{\mathcal{L}}),$$
which has a form similar to (\ref{regret}). The weight function on the right-hand side of the above equation is $g_1$, but as indicated below Remark~\ref{cororemark}, $\mathbb{E} \lambda_{g_r}(\mathcal{L}\;\Delta \; \wh{\mathcal{L}})$ has the same asymptotic form for all $r\geq 0$ up to a constant. The risk $\mathbb{E}\mathcal{R}(\wh{\mathcal{L}})$ has an ``empirical'' form 
\begin{align}\label{emprisk}
\wh{\mathcal{R}}_n(\wh{\mathcal{L}}) = \frac{1}{n} \sum_{i=1}^n e_{\wh{\mathcal{L}}}(X_i),
\end{align}
which can be used to evaluate the performance of a density level set estimator. Note that $e_{\wh{\mathcal{L}}}(X_i)$ still depends on the unknown $f$, unlike its counterparts for classification or regression level set (see Willett and Nowak, \cite{Willett07}). Nonetheless, we still use this ``empirical'' risk function as one of performance metrics in our simulation study, because the density functions are known there.
Our optimal bandwidth for level set estimation is based on the excess risk $\mathbb{E}\lambda_{g_r}(\mathcal{L}\;\Delta \; \wh{\mathcal{L}})$ for a $r\geq 0$, which is shown to resemble MISE for kernel density estimation in Theorem~\ref{riskexpect}, where we take $p=1$ and $g_r^{(p)}(x)=c^r\|\nabla f(x)\|$ (see Remark~\ref{assumpremark} b)(ii)). The following proposition provides another way of understanding this notion.
\begin{proposition}\label{L2limit}
When $g(x)=f(x)^r|f(x)-c|$ for some $r\geq 0$, under assumptions (K1), (K2), (F1), (H2), as $\delta\searrow0$ we have 
\begin{align}\label{L2AsLimit}
\frac{2\delta c^{-r}\mathbb{E}\lambda_g(\mathcal{L}\;\Delta \; \wh{\mathcal{L}})}{\mathbb{E} \int_{\mathcal{I}(\delta)} |\wh f(x) - f(x)|^2dx} \rightarrow 1.
\end{align} 
\end{proposition}

Following this result we can interpret the excess risk as a limit of the MISE for kernel density estimation constrained in a neighborhood of $\mathcal{M}$. As discussed in the remark after Corollary~\ref{L1asympexpress}, the risk $\mathbb{E}\lambda_{g}(\mathcal{L}\;\Delta \; \wh{\mathcal{L}})$ with $g\equiv 1$ or $g=f$ is analogous to MIAE used for kernel density estimation. In comparison with this $L^1$ type of risk, using the excess risk (which is $L^2$ type) for bandwidth selection in level set estimation enjoys some mathematical simplicity, similar to MISE for kernel density estimation (see page 16, Wand and Jones, \cite{Wand95}).  

In what follows we denote
\begin{align}
&m({\bf h}) = \mathbb{E} \int_{\cal{M}} \frac{|\wh f(x) - f(x)|^2}{\|\nabla f(x)\|}d\mathscr{H}(x), \label{MISEdef}\\
%
\text{and  } &\wt m({\bf h}) = s_n^2 \int_{\mathcal{M} } \frac{1}{\|\nabla f(x)\| } d\mathscr{H}(x) +  \int_{\mathcal{M} } \frac{\beta_{\mathbf{h}}(x)^2}{\|\nabla f(x)\| } d\mathscr{H}(x) . \label{MISEdef2}
\end{align}

The assumptions in Theorem~\ref{riskexpect} guarantee that when $g(x)=f(x)^r|f(x)-c|$ for some $r\geq 0$,
\begin{align*}
\mathbb{E}\lambda_g(\mathcal{L}\;\Delta \; \wh{\mathcal{L}}) = \frac{1}{2}c^rm({\bf h}) + o(\|{\bf h}\|^{2\nu} + s_n^2) = \frac{1}{2}c^r\wt m({\bf h}) + o(\|{\bf h}\|^{2\nu} + s_n^2).
\end{align*}
Therefore the excess risk can be asymptotically minimized by minimizing $\wt m({\bf h})$.
Note that
\begin{align}\label{AMISE}
\wt m({\bf h}) & =  \frac{1}{(\nu!)^2}\kappa_\nu^2 \sum_{k=1}^d \sum_{l=1}^d h_k^\nu h_l^\nu \int_{\cal{M}} \frac{f_{(k*\nu)}(x)f_{(l*\nu)}(x)}{\|\nabla f(x)\|}d\mathscr{H}(x)\nonumber\\
&\hspace{2cm} + \frac{1}{n\Pi_{j=1}^d h_j} \|K\|_2^2 \int_{\cal{M}}\frac{c}{\|\nabla f(x)\|}d\mathscr{H}(x)\nonumber\\
& =   \frac{1}{(\nu!)^2}\kappa_\nu^2({\bf h}^\nu)^TA(f){\bf h}^\nu + \frac{cb(f)\|K\|_2^2}{n} \frac{1}{(h_1^\nu h_2^\nu\cdots h_d^\nu)^{1/\nu}},
\end{align}
where we denote ${\bf h}^\nu = (h_1^\nu, h_2^\nu, \cdots, h_d^\nu)^T$,  $b(f) = \int_{\cal{M}}\|\nabla f(x)\|^{-1}d\mathscr{H}(x)$ and 
\begin{align*}
A(f) = [a_{kl}]_{1\leq k,l\leq d} \text{ with } a_{kl} =   \int_{\cal{M}} \frac{f_{(k*\nu)}(x)f_{(l*\nu)}(x)}{\|\nabla f(x)\|}d\mathscr{H}(x).
\end{align*}
For $\mathbf{u}=(u_1, u_2, \cdots, u_d)^T$, define the function
\begin{align}\label{Qfun}
Q(\mathbf{u};\mathbf{M},a,\nu) = \frac{1}{(\nu!)^2}\mathbf{u}^T\mathbf{M}\mathbf{u} + \frac{a}{(u_1u_2\cdots u_d)^{1/\nu}}.
\end{align}
Then from (\ref{AMISE}) we can write
\begin{align}\label{AMISEexpress}
\wt m({\bf h}) = Q \left({\bf h}^\nu;\kappa_\nu^2 A(f),\frac{cb(f)\|K\|_2^2}{n},\nu\right).
\end{align}

To ensure the uniqueness of the minimizer of $\wt m({\bf h})$, we impose the following assumption.
\begin{itemize}
\item[] {\bf (F2)} For $d=1$, we require $A(f)>0$; for $d\geq 2$, we require that $A(f)$ is positive semi-definite and $\inf_{\mathbf{u}\in\bar{\mathbb{R}}_+^d, \|\mathbf{u}\|\neq 0}\mathbf{u}^T A(f) \mathbf{u}/\|\mathbf{u}\|^2>0$, where $\bar{\mathbb{R}}_+^d$ denotes the set of vectors in $\mathbb{R}^d$ with non-negative coordinates.
\end{itemize}

An assumption similar to (F2) in the kernel regression setting appears in Yang and Tschernig (\cite{Yang99}). A density function with linearly dependent $\nu$th partial derivatives $\{f_{(k*\nu)}(x): k=1,\cdots,d\}$ does not satisfy this assumption. See Yang and Tschernig (\cite{Yang99}) for more discussions on the similar assumption.

Let $\mathbf{u}(\mathbf{M},a,\nu)$ be the vector $\mathbf{u}$ which minimizes $Q(\mathbf{u};\mathbf{M},a,\nu)$. Denote the $d\times d$ identity matrix by $\mathbf{I}_d$. We have the following optimization result for $\wt m({\bf h})$.
\begin{theorem}\label{optimAMISE}
Under assumptions (K1), (F1), and (F2), $\wt m({\bf h})$ is uniquely minimized by a bandwidth given by
\begin{align}\label{optimtheory}
\wt{ \mathbf{h}}_{\text{opt}} = \left( \frac{cb(f)\|K\|_2^2}{\kappa_\nu^2n} \right)^{1/(d+2\nu)} \mathbf{u}^{1/\nu}(A(f),1,\nu).
\end{align}
In addition, assume that $f$ has bounded and continuous $(\nu+2)$ times derivatives and $\int_{\mathbb{R}}|u^{\nu+2}\wt K(u)|du<\infty$. Then as $n\rightarrow\infty$, the bandwidth $\mathbf{h}_{\text{opt}}$ which minimizes $m({\bf h})$ satisfies
\begin{align}
\wt{ \mathbf{h}}_{\text{opt}} & =  \left\{\mathbf{I}_d+O\left(n^{-2\nu/(d+2\nu)}\right)\right\} \mathbf{h}_{\text{opt}},\label{hoptapprox} 
\end{align}
and
\begin{align}
\wt m(\wt{ \mathbf{h}}_{\text{opt}}) & =  \left\{1+O\left(n^{-2\nu/(d+2\nu)}\right)\right\} m(\mathbf{h}_{\text{opt}}).\label{mhoptapprox}
\end{align}
\end{theorem}
\begin{remark} \label{optimalbandform}$\;$\\[-14pt]
{\normalfont
\begin{itemize}
\item[a)] The result (\ref{optimtheory}) also contains the case $d=1$, which we state explicitly below. For $d=1$, we write
\begin{align}
& m(h) = \sum_{i=1}^N \frac{\mathbb{E}|\wh f(x_i) - f(x_i)|^2}{|f^\prime(x_i)|}, \\
 \text{ and } & \wt m(h) = \frac{\|K\|_2^2c}{nh}\sum_{i=1}^N \frac{1}{|f^\prime(x_i)|} + \frac{h^{2\nu}}{(\nu!)^2} \kappa_\nu^2 \sum_{i=1}^N \frac{f^{(\nu)}(x_i)^2}{|f^\prime(x_i)|},\label{asympriskd1}
\end{align}
where $f^{(\nu)}$ is the $\nu$th derivative of $f$. Then $m(h) = \wt m(h) + o(\frac{1}{nh} + h^4)$. The asymptotic optimal bandwidth is given by
\begin{align}\label{Constant}
\wt{h}_{\text{opt}} = Cn^{-\frac{1}{1+2\nu}} \text{ with } C = \left( \frac{c(\nu!)^2\|K\|_2^2 \sum_{i=1}^N |f^\prime(x_i)|^{-1} }{2\nu\kappa_\nu^2\sum_{i=1}^N [f^{(\nu)}(x_i)]^2|f^\prime(x_i)|^{-1}} \right)^{\frac{1}{1+2\nu}}.
\end{align}
\item[b)] If we have the restriction $h_1=h_2=\cdots h_d$ for $\mathbf{h}$, then $\wt{\mathbf{h}}_{\text{opt}}=(\wt h_{\text{opt}} , \cdots,\wt h_{\text{opt}} )^T$ has a closed form with
\begin{align*}
\wt h_{\text{opt}} = \left( \frac{cd (\nu!)^2 b(f) \|K\|_2^2}{ 2n\nu \kappa_\nu^2\sum_{k=1}^d\sum_{l=1}^d a_{kl} }\right)^{1/(d+2\nu)}.
\end{align*}
\item[c)] In general, for the multivariate case, (\ref{optimtheory}) has an analytical expression only when $d=2$, given by $\wt{ \mathbf{h}}_{\text{opt}} = (\wt h_{\text{opt},1} , \wt h_{\text{opt},2} )^T$, where
\begin{align*}
&\wt h_{\text{opt},1} = \left( \frac{c(\nu!)^2 b(f) \|K\|_2^2 a_{22}^{(\nu+1)/(2\nu)}}{2n\nu\kappa_\nu^2  a_{11}^{(\nu+1)/(2\nu)} ( a_{11}^{1/2} a_{22}^{1/2} + a_{12})} \right)^{1/(2+2\nu)} ,\\
 \text{and }\;& \wt h_{\text{opt},2} = \left(\frac{ a_{11}}{ a_{22}}\right)^{1/(2\nu)} \wt h_{\text{opt},1}.
\end{align*}
For $d\geq3$, one has to use numerical methods to find the solution. See Wand and Jones (\cite{Wand94}).
\end{itemize}
}
\end{remark}

Since (\ref{optimtheory}) contains unknown quantities, in practice we need to find estimators $\wh b(f)$ and $\wh A(f)$ for $b(f)$ and $A(f)$. Then the asymptotic risk function $\wt m(\mathbf{h})$ is estimated by
\begin{align}\label{mestimate}
\wh m({\bf h}) = Q \left({\bf h}^\nu;\kappa_\nu^2 \wh A(f),\frac{c\wh b(f)\|K\|_2^2}{n},\nu\right).
\end{align}
Correspondingly, the plug-in optimal bandwidth becomes
\begin{align}\label{optimprac}
\wh{\mathbf{h}}_{\text{opt}} = \left( \frac{c\wh b(f)\|K\|_2^2}{\kappa_\nu^2n} \right)^{1/(d+2\nu)} \mathbf{u}^{1/\nu}(\wh A(f),1,\nu).
\end{align}

For simplicity, below we assume $\nu=2$ in (\ref{AMISEexpress}) and (\ref{optimtheory}), but our methodology applies to general $\nu\geq 2$. Note that $b(f)$ and $A(f)$ involve the unknowns $\mathcal{M}$, $\nabla f$ and $f_{(k,k)}f_{(j,j)}$ for $1\leq k,j\leq d$, which need to be estimated. Below we discuss our choices of estimators and the relatively rates of convergence of our plug-in bandwidth selectors for $d=1$ and $d\geq 2$ separately, because the case $d\geq2$ involves estimation of surface integrals on level sets, whereas the case $d=1$ only requires point estimation.

We first consider $d=1$. Recall that $\mathcal{M} = \{x_i:\; i=1,2,\cdots, N\}$ for $d=1$ (see the discussion after the assumptions). Let $\wh{\mathcal{M}} = \{x: \wh f(x) = c\} = \{\wh x_i:\; i=1,2,\cdots, \wh N\}$, where $\wh N$ is the cardinality of $\wh{\mathcal{M}}$. Also let $\wh b(f) = \sum_{i=1}^{\wh N} |\wh f^\prime(\wh x_i)|^{-1}$ and $\wh A(f)=\sum_{i=1}^{\wh N} [\wh f^{\prime\prime}(\wh x_i)]^2|\wh f^\prime(\wh x_i)|^{-1}$. For $d=1$ and $\nu=2$, the estimated risk function in (\ref{mestimate}) is
\begin{align}\label{mestimate1}
\wh m({h}) = \frac{\|K\|_2^2c}{nh}\sum_{i=1}^{\wh N} \frac{1}{|\wh f^\prime(x_i)|} + \frac{h^{2\nu}}{(\nu!)^2} \kappa_\nu^2 \sum_{i=1}^{\wh N} \frac{[\wh f^{\prime\prime}(x_i)]^2}{|\wh f^\prime(x_i)|},
\end{align}
and the plug-in estimator in (\ref{optimprac})  is 
\begin{align}\label{optimprac1}
\wh{h}_{\text{opt}} = \wh C n^{-\frac{1}{5}} \text{ with } \wh C = \left( \frac{c\|K\|_2^2 \sum_{i=1}^{\wh N} |\wh f^\prime(\wh x_i)|^{-1} }{\kappa_2^2\sum_{i=1}^{\wh N} [\wh f^{\prime\prime}(\wh x_i)]^2|\wh f^\prime(\wh x_i)|^{-1}} \right)^{\frac{1}{5}}.
\end{align}

Note that in the above estimator we are essentially estimating $f$, $f^\prime$ and $f^{\prime}$ using $\wh f$, $\wh f^\prime$ and $\wh f^{\prime\prime}$. The kernel function $K$ may be replaced by a different one in these estimators. However, for simplicity of notation, we keep using $K$ in what follows. 
The bandwidths used in the estimators $\wh f$, $\wh f^\prime$ and $\wh f^{\prime\prime}$ can be chosen separately, for which we denote as ${h}^{(0)}$, ${h}^{(1)}$ and ${h}^{(2)}$, respectively. We propose to use the direct plug-in bandwidths for the kernel density and its first two derivatives as the pilot bandwidths ${h}^{(0)}$, ${h}^{(1)}$ and ${h}^{(2)}$, respectively. See Wand and Jones (\cite{Wand94}, \cite{Wand95}), Duong and Hazelton (\cite{Duong03}), and Chac\'{o}n et al. (\cite{Chacon11}) for details of the direct plug-in strategies. In fact, our pilot bandwidths for $d=1$ can be chosen following the exact procedure given in Samworth and Wand (\cite{Samworth10}, page 1777). The following theorem gives the relative rates of convergence of estimating our optimal bandwidth for $d=1$. Recall that $\wt m$ given in (\ref{asympriskd1}) is an asymptotic approximation to the excess risk when $g(x)=|f(x)-c|$ and $\wt{h}_{\text{opt}}$ given in (\ref{Constant}) is a minimizer of $\wt m$.

\begin{theorem}\label{practicald1}
Suppose $d=1$ and assumptions (F1), (F2), (K1) and (K2) hold with $\nu=2$. In addition, assume that $f$ has bounded continuous fourth derivatives and $K$ has bounded continuous third derivatives of bounded variation. If  $h^{(0)} \asymp n^{-1/5}$, $h^{(1)} \asymp n^{-1/7}$ and $h^{(2)} \asymp n^{-1/9}$, then for $\wh{h}_{opt}$ in (\ref{optimprac1}) and $\wh m$ in (\ref{mestimate1}) we have 
\begin{align}
&\wh{h}_{opt} = \wt{h}_{\text{opt}} \left\{1+O_p\left(n^{-2/9}\right)\right\},\label{hoptd1}\\
\text{and}\;\; & \wh m(\wh{h}_{opt})  = \wt m(\wt{h}_{\text{opt}}) \left\{1+O_p\left(n^{-2/9}\right)\right\}.\label{hoptd1m}
\end{align}
\end{theorem}
\begin{remark} $\;$\\[-14pt]
{\normalfont
\begin{itemize}
\item[] It is clear from the proof of the above theorem that the relative rates of convergence in (\ref{hoptd1}) and (\ref{hoptd1m}) are mainly determined by the choice of $h^{(2)}$ for the estimator $\wh f^{\prime\prime}$. If we choose $h^{(0)} =h^{(1)} =h^{(2)} \asymp n^{-1/9}$, then the conclusion in Theorem~\ref{practicald1} still holds.
\end{itemize}
}
\end{remark}

Next we consider $d\ge 2$. The asymptotics for $\wh{\mathbf{h}}_{\text{opt}}$ when $d\geq 2$ involves the estimation of integrals on level sets, which is studied in Qiao (\cite{Qiao19}). In the literature, estimating the volume of manifolds or surface integrals has been studied in, e.g., Cuevas et al. (\cite{Cuevas07}) and Jim\'{e}nez and Yukich (\cite{Jimenez11}). 
We consider plug-in estimators $\wh b(f) =  \int_{\wh{\cal{M}}}\|\nabla \wh f(x)\|^{-1}d\mathscr{H}(x)$ and $\wh A(f) = [\wh a_{kl}]_{1\leq k,l\leq d}$, where $\wh{\cal{M}} = \{x\in\mathbb{R}^d:\; \wh f(x)=c\}$ and 
\begin{align}\label{ConstantEst}
\wh a_{kl} = \int_{\wh{\cal{M}}} \|\nabla \wh f(x)\|^{-1}\wh f_{(k,k)}(x)\wh f_{(l,l)}(x)d\mathscr{H}(x).
\end{align}

Similar to the case $d=1$, we can still use different bandwidths for the estimation of derivatives of different orders in $\wh b(f)$ and $\wh A(f)$. Here for simplicity we choose to use a common bandwidth $\mathbf{h}_{\text{pilot}}$ in $\wh b(f)$ and $\wh A(f)$ for the reason given in the remark after Theorem~\ref{practicald1}. The following theorem is a consequence of Theorem 3.1 in Qiao (\cite{Qiao19}) by noticing that $\wt{\mathbf{h}}_{\text{opt}}$ and $\wt m(\wt{\mathbf{h}}_{\text{opt}})$ are smooth functions of $b(f)$ and $A(f)$, where $\wt m$ and $\wt{\mathbf{h}}_{\text{opt}}$ are given in (\ref{MISEdef2}) and (\ref{optimtheory}), respectively.

\begin{theorem}\label{practicald2}
Suppose $d\geq 2$ and assumptions (F1), (F2), (K1) and (K2) hold with $\nu=2$. In addition, assume that both $f$ and $K$ have continuous four times derivatives, and $K$ has bounded support. Let $h_n$ be a sequence such that $h_n\rightarrow0$ and $(\log{n})^{-1}nh_n^{d+4}\rightarrow\infty$ as $n\rightarrow\infty$. If $\mathbf{h}_{\text{pilot}} \asymp h_n$, then
\begin{align}
&\wh{\mathbf{h}}_{opt} = \wt{\mathbf{h}}_{\text{opt}} \left\{1+O_p\left(\alpha_{n}\right)\right\},\label{hoptd2}\\
\text{and}\;\; &\wh m(\wh{\mathbf{h}}_{opt})  = \wt m(\wt{\mathbf{h}}_{\text{opt}}) \left\{1+O_p\left(\alpha_{n}\right)\right\},\label{hoptd2m}
\end{align}
where $\alpha_{n} = \frac{1}{\sqrt{nh_n^5}} + \frac{1}{nh_n^{d+4}} + h_n^2$.
\end{theorem}

\begin{remark} \label{finalremark}$\;$\\[-14pt]
{\normalfont
\begin{itemize}
\item[a)] To minimize $\alpha_{n}$, we choose $\mathbf{h}_{\text{pilot}} \asymp n^{-1/\max\{9, d+6\})}$, i.e., $\mathbf{h}_{\text{pilot}} \asymp n^{-1/9}$ when $d=2$; and  $\mathbf{h}_{\text{pilot}} \asymp n^{-1/(d+6)}$ when $d\geq 3$. If so, then correspondingly we have $\alpha_{n}\asymp n^{-2/(\max\{9,d+6\})}$. In practice, we can use $\mathbf{h}^{(1)}$ as $\mathbf{h}_{\text{pilot}}$, which is the direct plug-in optimal bandwidth for estimating the gradient of $f$, because $\mathbf{h}^{(1)} \asymp n^{-1/(d+6)}$. If so, then we have $\alpha_{n}\asymp n^{-3/16}$ when $d=2$; and $\alpha_{n}\asymp n^{-2/(d+6)}$ when $d\geq 3$.
\item[b)] When $d\geq 3$, the computation of the surface integrals in $\wh b(f)$ and $\wh A(f)$ might be quite challenging. Alternatively, for a sequence $\epsilon_n>0$, we can replace $\wh b(f)$ and $\wh a_{kl}$ in $\wh{\mathbf{h}}_{\text{opt}}$ by the following two types of estimators using integration over small neighborhoods of $\wh{\cal{M}}$, where we still use $\mathbf{h}_{\text{pilot}}$ as the pilot bandwidth:  
\begin{align*}
&(i) 
\begin{cases}
\wh b^*(f) = \frac{1}{2\epsilon_n} \lambda(\wh f^{-1}[c-\epsilon_n,c+\epsilon_n])\\
\wh a_{kl}^* =\frac{1}{2\epsilon_n} \int_{\wh f^{-1}[c-\epsilon_n,c+\epsilon_n] } \wh f_{(k,k)}(x)\wh f_{(l,l)}(x) dx
\end{cases},\\
\text{  or }\\
&(ii) 
 \begin{cases}
 \wh b^\dagger(f) = \frac{1}{2\epsilon_n} \int_{\wh{\cal{M}}\oplus\epsilon_n}\|\nabla \wh f(x)\|^{-1}dx\\
 \wh a_{kl}^\dagger = \frac{1}{2\epsilon_n} \int_{\wh{\cal{M}}\oplus\epsilon_n} \|\nabla \wh f(x)\|^{-1}\wh f_{(k,k)}(x)\wh f_{(l,l)}(x)dx
 \end{cases}.
\end{align*}
Here $\epsilon_n$ controls the width of tubes around $\wh{\cal{M}}$ as domains of integration in these estimators. If we use these two types of estimators, then $\alpha_{n}$ is replaced by $\alpha_{n} + \epsilon_n^2$ in Theorem~\ref{practicald2} under the same condition. Again this is a consequence of Theorem 3.1 in Qiao (\cite{Qiao19}). Using $\epsilon_n$ of the same order of $h_n$ does not increase the previous relative rates of convergence $O_p(\alpha_{n})$. For example, if $\mathbf{h}_{\text{pilot}}$ is chosen to be $\mathbf{h}^{(1)}$ as in a), then we can use $\text{min}(\mathbf{h}^{(1)})$, $\text{max}(\mathbf{h}^{(1)})$, or the average of the individual bandwidths in $\mathbf{h}^{(1)}$ as $\epsilon_n.$ 
\end{itemize}
}
\end{remark}

\section{Simulation results}
A simulation study was run to assess the performance of our bandwidth selector $\wh{\mathbf{h}}_{\text{opt}}$ tailored for level set estimation. We compared the performance of our bandwidth selector with the least square cross validation method (see Rudemo, \cite{Rudemo82}, and Bowman, \cite{Bowman84}), which is an ISE-based selector denoted by $\mathbf{h}_{\text{LSCV}}$, as well as the direct plug-in bandwidth selector (see e.g. Wand and Jones, \cite{Wand94}) denoted by  $\mathbf{h}_{\text{DPI}}$. Note that both $\mathbf{h}_{\text{LSCV}}$ and $\mathbf{h}_{\text{DPI}}$ are bandwidth selectors for kernel density estimation. In order to make a fair comparison with $\wh{\mathbf{h}}_{\text{opt}}$, $\mathbf{h}_{\text{LSCV}}$ and $\mathbf{h}_{\text{DPI}}$ are also $d$-dimensional vectors, which correspond to diagonal bandwidth matrices.  

We first compared the performance of the three bandwidth selectors by considering a Gaussian mixture model with the distribution
\begin{align}\label{model13}
\frac{2}{3} \mathscr{N}\left( \begin{pmatrix} 0\\0\end{pmatrix}, \begin{pmatrix} 1/4 & 0\\0 & 1\end{pmatrix}\right) + \frac{1}{3} \mathscr{N}\left( \begin{pmatrix} 0\\0\end{pmatrix}, \frac{1}{50}\begin{pmatrix} 1/4 & 0\\0 & 1\end{pmatrix}\right),
\end{align}
which has a sharp mode and was constructed to represent a bivariate analog to density 4 in Marron and Wand (\cite{Marron92}). The levels of the density functions in our analysis were chosen corresponding to the 20\%, 50\% and 80\% HDRs, respectively, that is, $c = c(\tau)$, where $\tau=0.2, 0.5,$ and $0.8$ (see Section 1 for the definition of HDRs). 500 samples were drawn from this distribution, and for each sample we used the error $e(\mathbf{h}) = \lambda_{g_0}(\mathcal{L} \Delta \wh{\mathcal{L}})$, where $g_0(x)=|f(x)-c|$, to evaluate the performance of a density level set estimator with bandwidth $\mathbf{h}$ and the Gaussian kernel. 

\begin{figure}[!htbp]
\begin{center}
\includegraphics[scale=0.4]{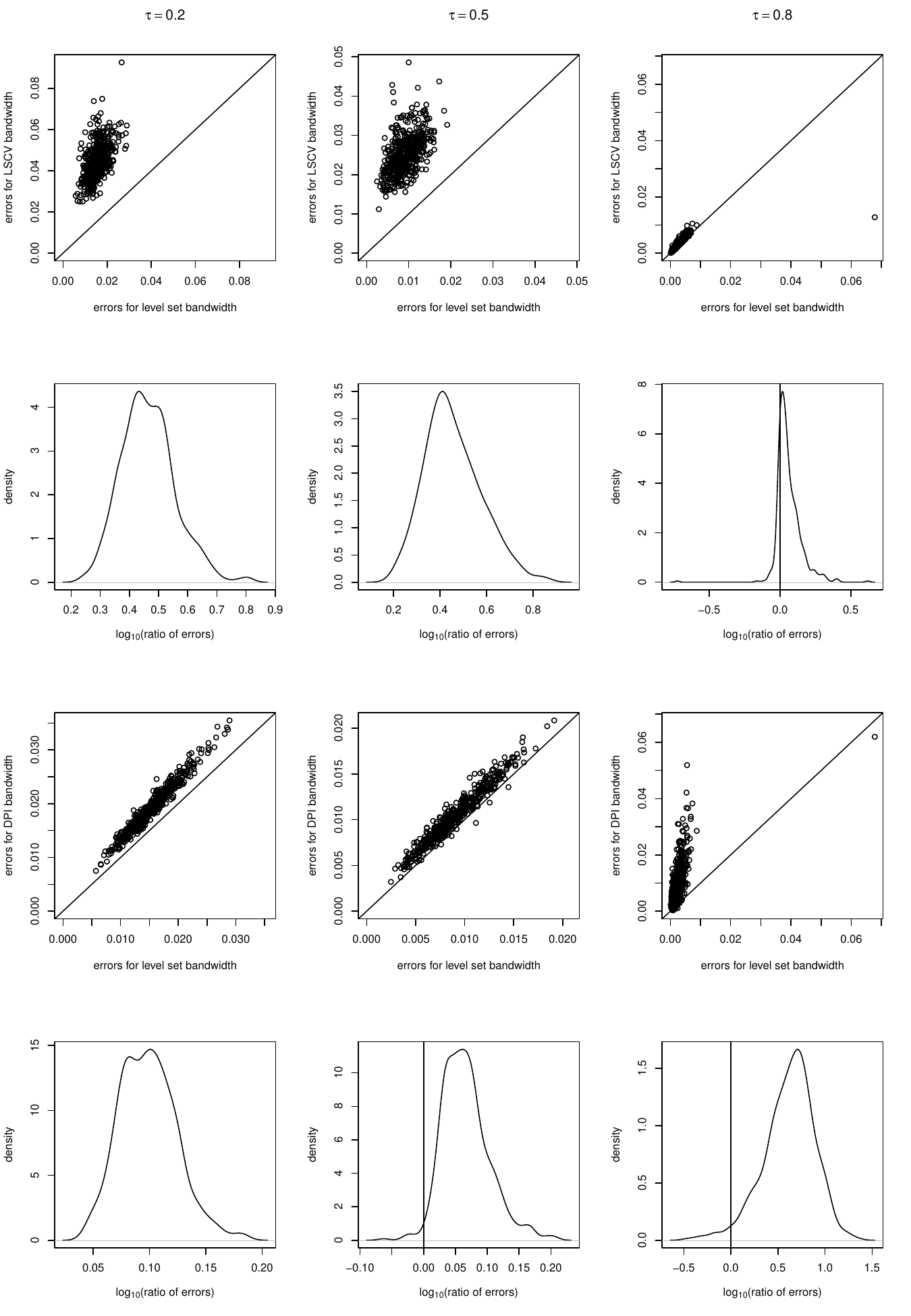}
\caption{Graphical comparison of the performance between $\wh{\mathbf{h}}_{\text{opt}}$ and $\mathbf{h}_{\text{LSCV}}$ and between $\wh{\mathbf{h}}_{\text{opt}}$ and $\mathbf{h}_{\text{DPI}}$ for $\tau=0.2, 0.5, 0.8$ for the model in (\ref{model13}), with sample size of $n=1000$ for 500 replications. The graphs in the first row show the scatter plots of the errors $e(\mathbf{h})$ for $\wh{\mathbf{h}}_{\text{opt}}$ and $\mathbf{h}_{\text{LSCV}}$, and the graphs in the second row show the kernel density estimates of the common logarithm of ratios between the errors using $\mathbf{h}_{\text{LSCV}}$ and the errors using $\wh{\mathbf{h}}_{\text{opt}}$. The third and fourth rows are similar comparisons between $\wh{\mathbf{h}}_{\text{opt}}$ and $\mathbf{h}_{\text{DPI}}$.}
\label{fig:simratios1}
\end{center}
\end{figure}

Figure \ref{fig:simratios1} show the simulation results for the model in (\ref{model13}) with sample size $n=1000$. It can be seen that our bandwidth selector $\wh{\mathbf{h}}_{\text{opt}}$ performed better than $\mathbf{h}_{\text{LSCV}}$ and $\mathbf{h}_{\text{DPI}}$ in terms of the error $e(\mathbf{h})$ for most of the samples.  The improvement of $\wh{\mathbf{h}}_{\text{opt}}$ was statistically significant at the $0.1\%$ level for $\tau=0.2, 0.5$, and $0.8$ when the Wilcoxon tests were applied to the ratio of errors given by $e(\mathbf{h}_{\text{LSCV}})/e(\wh{\mathbf{h}}_{\text{opt}})$ and $e(\mathbf{h}_{\text{DPI}})/e(\wh{\mathbf{h}}_{\text{opt}})$. For each $\tau$ value, among the 500 samples, we chose the one with the ratio of errors $e(\mathbf{h}_{\text{DPI}})/e(\wh{\mathbf{h}}_{\text{opt}})$ closest to the median as a representative. 

\begin{figure}[h]
\begin{center}
\includegraphics[scale=0.46]{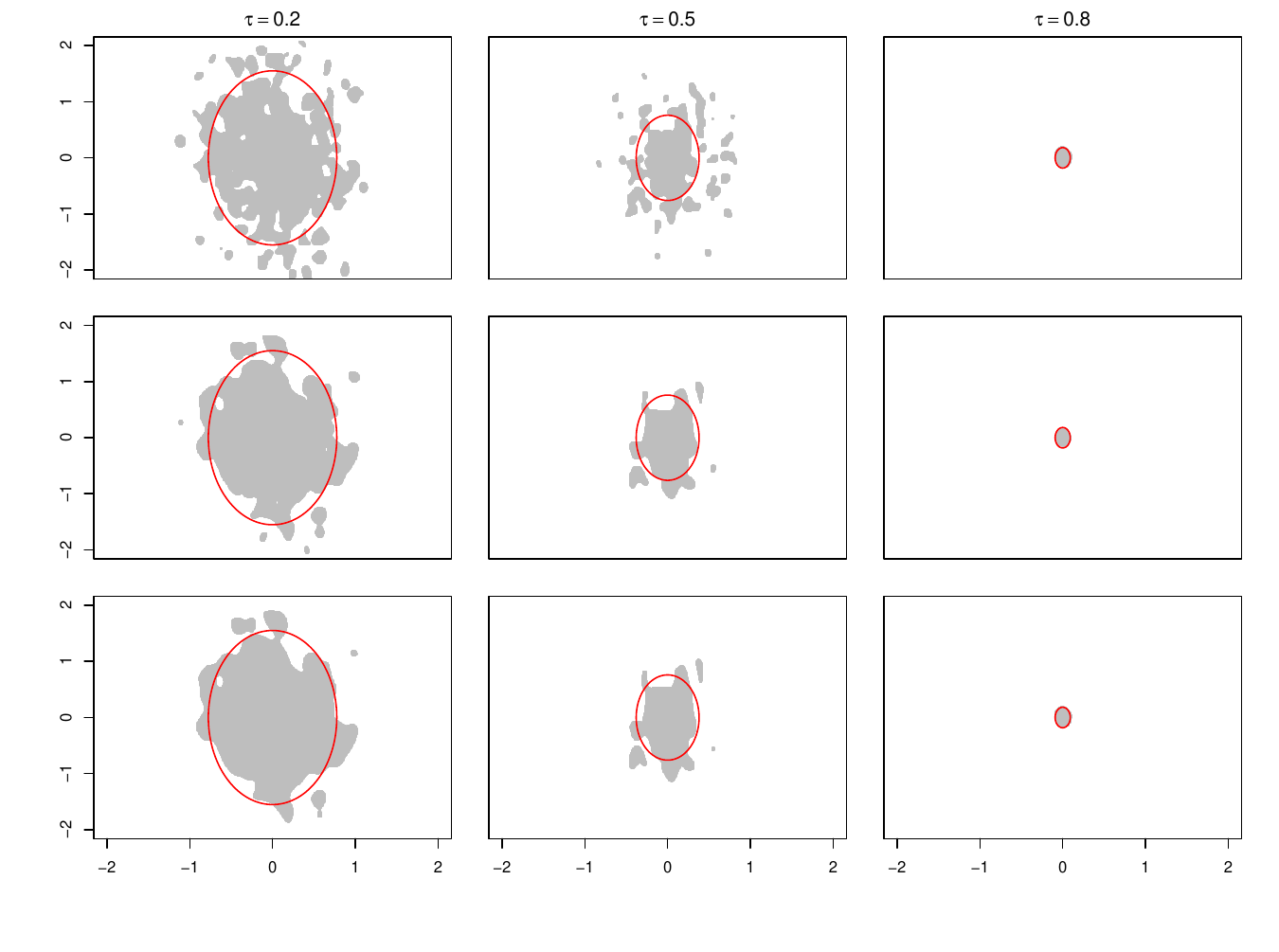}
\caption{Comparisons among the estimated level sets using $\mathbf{h}_{\text{LSCV}}$ (upper panels), $\mathbf{h}_{\text{DPI}}$ (middle panels) and $\wh{\mathbf{h}}_{\text{opt}}$ (lower panels) for $\tau=0.2, 0.5, 0.8$ for the model in (\ref{model13}). Estimated level sets are represented by the gray areas, and the true level sets are enclosed by the red curves. The samples were chosen such that the ratios between the errors using $\mathbf{h}_{\text{DPI}}$ and the errors using $\wh{\mathbf{h}}_{\text{opt}}$ are closest to their medians in the 500 replications.}
\label{fig:estimation}
\end{center}
\end{figure}

Figure \ref{fig:estimation} visually compares the level set estimations between the three bandwidth selectors for the representative samples for $\tau=0.2, 0.5$, and $0.8$. It can be seen that when $\wh{\mathbf{h}}_{\text{opt}}$ or $\mathbf{h}_{\text{DPI}}$ were used, the level sets were estimated reasonably well, with $\wh{\mathbf{h}}_{\text{opt}}$ slightly better, while using $\mathbf{h}_{\text{LSCV}}$ only captured the level sets for $\tau=0.8$. When we decreased the sample size to $n=500$, $\wh{\mathbf{h}}_{\text{opt}}$ still performed better than $\mathbf{h}_{\text{LSCV}}$ for $\tau=0.2$, $0.5$ and $0.8$, and better than $\mathbf{h}_{\text{DPI}}$ for $\tau=0.2$ and $0.8$ but not for $\tau = 0.5$. With this smaller sample size $\wh{\mathcal{M}}$ using the pilot bandwidth had about 17\% chance to be an empty set in the replications for $\tau=0.8$, which corresponds to a relatively high density level, and in these cases our bandwidth selector was not computable and so we had set $e(\mathbf{h}) = \lambda_g(\mathcal{L})$. This issue arises because the kernel density estimator underestimates the density in a neighborhood of the modes on average, when a second order kernel is used (see the expansion of the bias in (\ref{MISEbias})). When the level $c$ is relatively high and the sample size $n$ is small, we suspect that using a higher order kernel or a more sophisticated pilot bandwidth in the pilot density estimate might make an improvement on this issue. 
%

In addition, we also considered 12 bivariate Gaussian mixture models used in Wand and Jones (\cite{Wand93}), which cover from unimodal to quadrimodal models. We further extended these density functions to their trivariate counterparts in our simulation study as specified below. Denote a bivariate Gaussian mixture model with $k\geq 1$ components by $\sum_{i=1}^k w_i \mathscr{N}(\mu_i, \Sigma_i)$, where for some $-1<\rho_i<1,$
\begin{align*}
\mu_i = \begin{pmatrix} \mu_{i1}\\ \mu_{i2} \end{pmatrix}, \text{ and } \Sigma_i = \begin{pmatrix}
\sigma_{i1}^2 & \rho_i\sigma_{i1} \sigma_{i2} \\
 \rho_i\sigma_{i1} \sigma_{i2}  & \sigma_{i2}^2
\end{pmatrix}.
\end{align*} 
Its trivariate extension is $\sum_{i=1}^k w_i \mathscr{N}(\wt\mu_i, \wt\Sigma_i)$, where
\begin{align*}
\wt\mu_i = \begin{pmatrix} \mu_{i1}\\ \mu_{i2} \\ \mu_{i2} \end{pmatrix}, \text{ and } \wt\Sigma_i = \begin{pmatrix}
\sigma_{i1}^2 & \rho_i\sigma_{i1} \sigma_{i2} & \rho_i\sigma_{i1} \sigma_{i2} \\
 \rho_i\sigma_{i1} \sigma_{i2}  & \sigma_{i2}^2 & \rho_i\sigma_{i2} \sigma_{i2}\\
 \rho_i\sigma_{i1} \sigma_{i2} & \rho_i\sigma_{i2} \sigma_{i2} & \sigma_{i2}^2
\end{pmatrix}.
\end{align*} 
In other words, the third marginal means and variances replicate the second ones for all the components and the correlation coefficients remain the same. If $\wt\Sigma_i$ is not positive definite by this extension, we replaced $\rho_i$ by its half, which makes the covariance matrices of all the components positive definite. 

These 12 models and their extensions are used to compare the performance of the three bandwidth selectors in the following four cases, where case 1 can be viewed as a base case, and in cases 2-4 we consider the variation of the risk criteria, the dimensions $d$, and the orders of the kernel function $\nu$, respectively. Recall that $e(\mathbf{h})$ denotes a error metric for density level set estimation with bandwidth $\mathbf{h}$.  

Case 1: d=2, $\nu=2$, using $e(\mathbf{h}) = \lambda_{g_0}(\mathcal{L} \Delta \wh{\mathcal{L}})$;

Case 2: d=2, $\nu=2$, using $e(\mathbf{h}) = \wh{\mathcal{R}}_n(\wh{\mathcal{L}})$ as defined in (\ref{emprisk});

Case 3: d=2, $\nu=4$, using $e(\mathbf{h}) = \lambda_{g_0}(\mathcal{L} \Delta \wh{\mathcal{L}})$;

Case 4: d=3, $\nu=2$, using $e(\mathbf{h}) = \lambda_{g_0}(\mathcal{L} \Delta \wh{\mathcal{L}})$.\\[-10pt]

Before we show the simulation results, we give some details in the implementation. The Gaussian kernel was used for cases 1,2, and 4 (i.e., $\nu=2$). For case 3, the fourth-order kernel function is chosen to be $K(x_1,x_2)=\wt K(x_1)\wt K(x_2)$ with $\wt K(v) = \frac{1}{2}(3-v^2)\phi(v),\;v\in\mathbb{R},$ where $\phi$ is the pdf of a standard normal distribution. We used $\mathbf{h}^{(1)}$ as the pilot bandwidth $\mathbf{h}_{\text{pilot}}$ for $\nu=2$, which has been discussed in Remark~\ref{finalremark} a). For $\nu=4$, while we still use $\mathbf{h}^{(1)}$ to estimate $\wh f$ and $\nabla \wh f$, we use $\mathbf{h}^{(2)}$ as the pilot bandwidth to estimate the fourth derivatives of $f$, where $\mathbf{h}^{(2)}$ is the direct plug-in optimal bandwidth for estimating the Hessian of $f$. Our bandwidth selector $\wh{\mathbf{h}}_{\text{opt}}$ involves the calculation of line/surface integrals. The numerical approximation to line integrals on curves when $d=2$ are straightforward (for cases 1, 2, and 3). For case 4, we generated meshes with fine triangulation and used the corresponding Riemann sums to approximate the surface integrals. As indicated in Remark~\ref{finalremark} b), these surface integrals can also be approximated by integration over some small neighborhoods of the surfaces. 

For each of the distributions, random sampling was replicated for 500 times. Again we used $c=c(\tau)$ with $\tau=0.2, 0.5,$ and $0.8$ as the levels of the density functions. For each case, we have 36 combinations of the $\tau$ values and models. The sample sizes were chosen to be $n=1,000$, $n=2,000$, and $n=10,000$. We applied the one-sided Wilcoxon tests to the ratios of errors given by $e(\mathbf{h}_{\text{LSCV}})/e(\wh{\mathbf{h}}_{\text{opt}})$ and $e(\mathbf{h}_{\text{DPI}})/e(\wh{\mathbf{h}}_{\text{opt}})$, respectively. Table~\ref{simulationres} below summarizes the counts of scenarios when the improvement of $\wh{\mathbf{h}}_{\text{opt}}$ was not statistically significant at the $0.1\%$ levels for $\tau=0.2, 0.5$, and $0.8$. 

\begin{table}[h]
\centering
\caption{Simulation results}
\label{simulationres}
\resizebox{4.7in}{!}{
\begin{tabular}{ll|llll|llll}
\hline\hline
                        &      & \multicolumn{4}{c|}{$\wh{\mathbf{h}}_{\text{opt}}$ vs. $\mathbf{h}_{\text{LSCV}}$} & \multicolumn{4}{c}{$\wh{\mathbf{h}}_{\text{opt}}$ vs. $\mathbf{h}_{\text{DPI}}$} \\
                        &      & $\tau=0.2$  & $\tau=0.5$ & $\tau=0.8$ & Total  & $\tau=0.2$ & $\tau=0.5$ & $\tau=0.8$ & Total \\ \hline\hline
\multirow{3}{*}{case 1} & n=1000   & 4    & 3   & 2   & 9    & 10  & 8   & 2  & 20     \\
                        & n=2000 & 3    & 3   & 2  & 8       & 9   & 8   & 0  & 17     \\
                        & n=10000  &3     &1    &1   &5       &6    &2    &1   &9     \\ \hline\hline
\multirow{3}{*}{case 2} & n=1000 & 3    & 3   & 1   & 7      & 10  & 8   & 1  & 19     \\
                        & n=2000 & 2    & 3   & 1  & 6       & 9   & 8   & 1 & 18      \\
                        & n=10000  &3     &1    &1   &5       &7    &3    &0   &10      \\ \hline\hline
\multirow{3}{*}{case 3} & n=1000 & 6    & 4   & 3 & 13        & 7   & 5   & 1  & 13    \\
                        & n=2000   & 4    & 3   & 2   & 9     & 7   & 3   & 0  & 10    \\
                        & n=10000 & 3     & 1   & 0   &4       &3    &1    &0   &4      \\ \hline\hline
\multirow{3}{*}{case 4} & n=1000  & 11   & 10  & 4 & 25        & 11  & 11  & 2 & 24    \\
                        & n=2000 & 6    & 5   & 1   & 12      & 7   & 7   & 0  & 14    \\
                        & n=10000  &3     &3    &1   &7       &7    &3    &1   &11      \\  \hline\hline
\end{tabular}
}
\end{table}

Overall we find our bandwidth selector $\wh{\mathbf{h}}_{\text{opt}}$ performs better than $\mathbf{h}_{\text{LSCV}}$ and $\mathbf{h}_{\text{DPI}}$ for density level set estimation, especially when the sample size is moderately large. Between the two competitors $\mathbf{h}_{\text{LSCV}}$ and $\mathbf{h}_{\text{DPI}}$, in general a larger sample size is needed for $\wh{\mathbf{h}}_{\text{opt}}$ to outperform the latter, though the needed sample size can be reduced by using higher order kernels as shown in case 3. Also it appears $\wh{\mathbf{h}}_{\text{opt}}$ performs well for high density levels, while we need larger sample sizes for the asymptotics for $\wh{\mathbf{h}}_{\text{opt}}$ to show effect when the levels are low. Note that data of large sample sizes are available for many application areas of density level set estimation, such as flow cytometry (Naumann and Wand, \cite{Naumann09}) and astronomical survey (Jang, \cite{Jang06}). 

\section{Discussion}
In this paper we give asymptotic $L^p$ approximations of $\lambda_g (\mathcal{L}\;\Delta \; \wh{\mathcal{L}})$ and $\mathbb{E}\lambda_g (\mathcal{L}\;\Delta \; \wh{\mathcal{L}})$, where $p$ is determined by the local behavior of $g$ around $\mathcal{M}$. In particular, when $g(x) = f(x)^r|f(x)-c|$ for some $r\geq0$, the excess risk $\mathbb{E}\lambda_g (\mathcal{L}\;\Delta \; \wh{\mathcal{L}})$ has an $L^2$ approximation and is used to select bandwidth for density level set estimation. Numerical results verify that our bandwidth selectors tailored for level set estimation outperforms the lease square cross validation and the direct plug-in bandwidth selectors for density estimation, when the sample size is moderately large. 

As indicated in the Introduction section, the work in Doss and Weng (\cite{Doss18}) is related to some of the results in this paper, and they have given a comparison between their work with an earlier arXiv version of this paper. When focusing on the level set estimation, they only consider $g=f$ in the asymptotic approximation for $\mathbb{E} \lambda_g (\mathcal{L}\;\Delta \; \wh{\mathcal{L}})$ (which is an $L^1$ type of risk as interpreted in our Corollary~\ref{L1asympexpress}), and use it as a risk function for bandwidth selection for density level set estimation. We give the expressions of the asymptotic forms of both $\lambda_g (\mathcal{L}\;\Delta \; \wh{\mathcal{L}})$ and $\mathbb{E}\lambda_g (\mathcal{L}\;\Delta \; \wh{\mathcal{L}})$ for a general class of $g$, which allows us to interpret them as asymptotic $L^p$ type of loss and risk, depending on a property of $g$ given in assumption (G1). In these approximations higher order kernel functions are also allowed if higher smoothness of the density function is assumed. Our bandwidth selection is based on an $L^2$ type of risk (the excess risk), which corresponds to a specific choice of $g$ in our general result. The excess risk resembles the MISE for kernel density estimation, and is more tractable than the $L^1$ type of risk. Note that in order to study the theory for the  minimization of the $L^1$ type of risk (when $g=f$), Doss and Weng (\cite{Doss18}) assume that the density function $f$ is unimodal and symmetric (see their Corollary 2.1). By contrast, the minimization of the excess risk does not require such assumptions on the shape of the density functions, and our optimal bandwidth can have analytical forms, depending on the structure of the bandwidth matrix (see Remark~\ref{optimalbandform}). Also as indicated in Remark~\ref{cororemark}, if one uses  the $L_1$ type of risk (e.g. $g=f$ or $g\equiv1$) for bandwidth selection, minimizing one of its upper bounds is an alternative approach, which can give closed-form solutions (again depending on the structure of the bandwidth matrix). 

Note that $g(x) = f(x)^r|f(x)-c|$ for $r\geq 0$ is chosen to be used mainly because its close connection to the excess risk in the classification literature and its interpretation as local MISE given in Proposition~\ref{L2limit}. One can also use $g(x)=\|\nabla f(x)\| \,f(x)^r|f(x)-c|$, which adds the norm of the gradient as a weight into the integrand. If so then we have a simple approximation $\mathbb{E}\lambda_g(\mathcal{L}\;\Delta \; \wh{\mathcal{L}})\approx \frac{1}{2} c^r\mathbb{E} \int_{\cal{M}} |\wh f(x) - f(x)|^2d\mathscr{H}(x)$ by Theorem~\ref{riskexpect}. Using this risk, the asymptotic optimal bandwidth has a simpler form because $b(f)$ and $a_{kl}$ in $\wt m(\mathbf{h})$ in (\ref{AMISE}) will be replaced by $\mathscr{H}(\mathcal{M})$ and $\int_{\mathcal{M}}f_{(k*\nu)}(x)f_{(l*\nu)}(x)d\mathscr{H}(x)$, respectively. In other words, one does not need to estimate the first derivatives in the plug-in bandwidth selector. With this form of the surface integrals, one may use U-statistic type estimators to improve the relative rates of convergence given in Theorem~\ref{practicald2}. See Theorem 2.1 and Corollary 3.1 in Qiao (\cite{Qiao19}).

In this paper we focus on the selection of bandwidth vectors, which correspond to diagonal bandwidth matrices. We expect that our results can be extended to full unconstrained bandwidth matrices, which might work better for level set estimation. See, e.g. Chac\'{o}n and Duong (\cite{Chacon10}). A closely related question is to select bandwidths for the estimation of HDR. Here the risk criterion can be set as $\mathbb{E} \lambda_g (\mathcal{L}_{c(\tau)}\;\Delta \; \wh{\mathcal{L}}_{\wh c(\tau)})$ (see the Introduction section). Doss and Weng (\cite{Doss18}) use $g=f$ and thus obtain an $L_1$ approximation of this risk criterion. It is expected that one can have an $L_2$ approximation to this risk using $g(x) = f(x)^r|f(x)-c|$ for $r\geq 0$ and select bandwidths using similar ideas in this paper. We leave the exploration of this idea to future work.

\section{Proofs}

\begin{proof}[Proof of Theorem \ref{L1risk}]$\;$

We first present the proof for the case $d\geq2$. The case $d=1$ is briefly discussed at the end. 

By Theorem 2 of Cuevas et al. (\cite{Cuevas06}), we have
\begin{align}\label{hausdist}
d_H(\mathcal{M}, \widehat{\mathcal{M}}) = O\left( \|\wh f -f\|_\infty\right).
\end{align}

With a slight generalization of Theorem 1 of Einmahl and Mason (\cite{Einmahl05}) to the case of individual bandwidth for each dimension in the kernel density estimator, it follows from assumptions (K1), (K2), (F1) and (H1) that
\begin{align}\label{stochrate}
\limsup_{n\rightarrow\infty} \sqrt{\frac{nh_1\cdots h_d}{\log{n}}} \sup_{x\in\mathbb{R}^d}|\wh f(x) - \mathbb{E} \wh f(x)| \leq \eta_1,\; a.s.
\end{align}
for some constant $\eta_1$. Following a standard derivation for kernel density estimation, we can show that, there exists a constant $\eta_2>0$ such that 
\begin{align}\label{biasuniform}
\sup_{x\in\mathbb{R}^d} |\mathbb{E} \wh f(x) - f(x)| \leq \eta_2\|\mathbf{h}\|^\nu.
\end{align}

Combining (\ref{hausdist}), (\ref{stochrate}) and (\ref{biasuniform}), we have that there exists a positive constant $\eta_3$ such that with
\begin{align}\label{epsilon_n}
\epsilon_n = \eta_3\left(\sqrt{ \frac{\log n}{nh_1\cdots h_d}} + \|\mathbf{h}\|^\nu\right),
\end{align}
we have $\mathcal{L}\;\Delta \; \wh{\mathcal{L}} \subset \mathcal{M}\oplus\epsilon_n$ for $n$ large enough with probability one, and as a result,
\begin{align}\label{reducedset}
\lambda_g (\mathcal{L}\;\Delta \; \wh{\mathcal{L}}) = \int_{\mathbb{R}^d}  \mathbf{I}(x\in \mathcal{L}\;\Delta \; \wh{\mathcal{L}}) g(x) dx =\int_{\mathcal{M}\oplus\epsilon_n}  \mathbf{I}(x\in \mathcal{L}\;\Delta \; \wh{\mathcal{L}}) g(x) dx.
\end{align}
By the definition of reach for a manifold, for any $\epsilon$ with $0<\epsilon<\rho(\mathcal{M})$, we can write $\mathcal{M}\oplus\epsilon = \left\{\zeta_x(s):\; x\in\mathcal{M}, \; |s|\leq \epsilon\right\}$. Then for large $n$, the map $\zeta(x,s):=\zeta_x(s)$ is a diffeomorphism from $\mathcal{M}\times[-\epsilon_n,\epsilon_n]$ to $\mathcal{M}\oplus\epsilon_n$.

The following derivation uses integration on manifolds, the theory of which can be found in e.g., Guillemin and Pollack (\cite{Guillemin74}), page 168 and Gray (\cite{Gray04}), Theorems 3.15 and 4.7. Denote $\mathcal{S}_n=\mathcal{M}\times[-\epsilon_n,\epsilon_n]$. For any $(x,s)\in \mathcal{S}_n$, let $\psi: U\mapsto \mathcal{S}_n$ be a local parameterisation of $\mathcal{S}_n$ around $(x,s)$, where $U$ is an open subset of $\mathbb{R}^d$. Denote function composition $\eta=\zeta\circ\psi$, which is a local parameterisation of $\mathcal{M}\oplus\epsilon_n$ around $\zeta_x(s)$. Note that both $\psi$ and $\eta$ depend on $(x,s)$ and this dependence has been suppressed in our notation. Let $D\psi$ and $D\eta$ be the Jacobian matrices of $\psi$ and $\eta$, respectively. Define the derivative $D\zeta_{(x,s)}:T_{(x,s)}\mathcal{S}_n\mapsto \mathbb{R}$ of $\zeta$ at $(x,s)$ by $D\zeta_{(x,s)}=D\eta\circ(D\psi)^{-1}$, where $T_{(x,s)}\mathcal{S}_n=T_x\mathcal{M}\times\mathbb{R}$ is the tangent space of $\mathcal{S}_n$ at $(x,s)$. Following Proposition 6 in Cannings et al. (\cite{Cannings17}), we have
\begin{align}\label{jacobian}
D\zeta_{(x,s)}(v_1,v_2) = (I + s B(x)) \left(v_1+\frac{\nabla f(x)}{\|\nabla f(x)\|}v_2\right),
\end{align}
for $v_1\in T_x(\mathcal{M})$ and $v_2\in\mathbb{R}$, where $B(x) = \frac{1}{\|\nabla f(x)\|} \left( I - \frac{\nabla f(x) \nabla f(x)^T}{\|\nabla f(x)\|^2} \right) \nabla^2 f(x)$ is also known as the shape operator on level sets (see Qiao, \cite{Qiao19}). It then follows from (\ref{reducedset}), (\ref{jacobian}), and the derivation in Section 7.3 of Cannings et al. (\cite{Cannings17}) that with probability one for $n$ large enough we have
%
%
\begin{align}\label{detIbx}
\lambda_g (\mathcal{L}\;\Delta \; \wh{\mathcal{L}}) =  \int_{\mathcal{M}} \int_{-\epsilon_n}^{\epsilon_n} \text{det}(I + s B(x)) \; g(\zeta_x(s)) \; \mathbf{I}(\zeta_x(s)\in \mathcal{L}\;\Delta \; \wh{\mathcal{L}})  ds d\mathscr{H}(x).
\end{align}
%
Since $\text{det}(I + s B(x)) = 1 + o(1)$ as $n\rightarrow\infty$, uniformly in $s\in[-\epsilon_n,\epsilon_n]$ and $x\in\mathcal{M}$, we obtain
\begin{align}\label{reducedset2}
\lambda_g (\mathcal{L}\;\Delta \; \wh{\mathcal{L}}) =  \int_{\mathcal{M}} \int_{-\epsilon_n}^{\epsilon_n} g(\zeta_x(s)) \; \mathbf{I}(\zeta_x(s)\in \mathcal{L}\;\Delta \; \wh{\mathcal{L}})  ds d\mathscr{H}(x) \{1+o(1)\}.
\end{align}

For any $x\in\mathcal{M}$, recall that $P_n(x) = \zeta_x(t_n(x))\in \wh{\mathcal{M}}$. Using Lemma 1 in Chen et al. (\cite{Chen17}), $\mathcal{M}$ and $\wh{\mathcal{M}}$ are normal compatible and hence $P_n$ is well defined. For $n$ large enough, we have $\text{sign}(f(\zeta_x(s))-c)=\text{sign}(s)$ and $\text{sign}(\wh f(\zeta_x(s))-c)=\text{sign}(s-t_n(x))$, for $s\in(-\epsilon_n,\epsilon_n)$. Hence for any $x\in\mathcal{M}$, the event $\zeta_x(s)\in\mathcal{L}\;\Delta \; \wh{\mathcal{L}}$ is equivalent to $s\in[t_n(x)\wedge0, t_n(x) \vee 0]$, where $t_n(x)\wedge0=\min(t_n(x),0)$ and $t_n(x)\vee0=\max(t_n(x),0)$. With probability one we have $t_n(x)<\epsilon_n$ for $n$ large enough since $\mathcal{L}\;\Delta \; \wh{\mathcal{L}} \subset \mathcal{M}\oplus\epsilon_n$ as indicated above. Hence from (\ref{reducedset2}) we can write 
\begin{align}\label{lambdagapp}
\lambda_g (\mathcal{L}\;\Delta \; \wh{\mathcal{L}}) & = \int_{\mathcal{M}} \int_{-\epsilon_n}^{\epsilon_n}   \mathbf{I} (s\in[t_n(x)\wedge0, t_n(x) \vee 0]) \;g\left(\zeta_x(s) \right) ds  d\mathscr{H}(x)\{1+o(1)\}\nonumber\\
& =  \int_{\mathcal{M}} \int_{t_n(x)\wedge0}^{t_n(x) \vee 0} g\left(\zeta_x(s) \right) ds  d\mathscr{H}(x)\{1+o(1)\}.
\end{align}
By Assumption (G1) we have
\begin{align}\label{gapprox}
g\left(\zeta_x(s) \right) = g^{(p)}(x)|s|^p \{1+o(1)\},
\end{align}
where $o(1)$ is uniform in $(x,s)$ for $x\in\mathcal{M}$ and $s\in[-\epsilon_n,\epsilon_n]$. Note that
\begin{align}\label{integralapprox}
\int_{\mathcal{M}} \int_{t_n(x)\wedge0}^{t_n(x) \vee 0}  g^{(p)}(x)|s|^p ds  d\mathscr{H}(x)  = & \frac{1}{p+1}\int_{\mathcal{M}}g^{(p)}(x) |t_n(x)|^{p+1}  d\mathscr{H}(x). 
\end{align}
By using the Taylor expansion for $x\in\mathcal{M}$, we have
\begin{align*}
0=\wh f(P_n(x)) - f(x) = &\wh f(x) -f(x) + \frac{\nabla f(x)^T\nabla \wh f(x)}{\|\nabla f(x)\| } t_n(x) \\
&\hspace{1.8cm }+ O\left(|t_n(x)|^2\sup_{x\in\mathcal{M}\oplus\epsilon_n}\|\nabla^2\wh f(x)\|\right).
\end{align*}
It follows that  
\begin{align}\label{distapprox}
|t_n(x)| = \frac{|\wh f(x) - f(x)|}{\|\nabla f(x)\|} \{1+o_p(1)\},
\end{align}
where $o_p(1)$ is uniform in $x\in\mathcal{M}$. Also see Lemma 2.2 in Qiao (\cite{Qiao19}). Combining (\ref{integralapprox}) and (\ref{distapprox}), we have
\begin{align}\label{integralapprox2}
&\int_{\mathcal{M}} \int_{t_n(x)\wedge0}^{t_n(x) \vee 0}  g^{(p)}(x)|s|^p ds  d\mathscr{H}(x) \nonumber\\
= & \frac{1}{p+1}\int_{\mathcal{M}}\frac{g^{(p)}(x)}{\|\nabla f(x)\|^{p+1}} |\wh f(x) -f(x)|^{p+1}  d\mathscr{H}(x) \{1+o_p(1)\}.
\end{align}
%
The result (\ref{L1expexpress}) immediately follows from (\ref{lambdagapp}), (\ref{gapprox}) and (\ref{integralapprox2}). 

The proof for the case $d=1$ can be shown by going through a similar procedure as above, but should be simplified with fewer geometric ingredients. Note that (\ref{reducedset}) is still valid for $d=1$. Then using (\ref{gapprox}) we have
\begin{align*}
 \lambda_g (\mathcal{L}\;\Delta \; \wh{\mathcal{L}}) 
= & \sum_{x\in\mathcal{M}} \int_{t_n(x)\wedge0}^{t_n(x) \vee 0} g\left(\zeta_x(s) \right) ds \\
= & \frac{1}{p+1}\sum_{x\in\mathcal{M}} \frac{g^{(p)}(x)}{|f^\prime(x)|^{p+1}} |\wh f(x) -f(x)|^{p+1}\{1+o_p(1)\}. 
\end{align*}
$\hfill\square$
\end{proof}


\begin{proof}[Proof of Theorem \ref{riskexpect}]$\;$

We only show the proof for the case $d\geq 2$, as the proof is similar and simpler for the case $d=1$, as shown in the proof of Theorem \ref{L1risk}. Before we show the main steps in the proof, we give a useful property of the kernel function $K$ under assumption (K1):
\begin{align}\label{Kproperty}
K\in\mathscr{L}^q, \text{ for all } 1\leq q\leq \infty.
\end{align}
%
To show (\ref{Kproperty}), notice that for $1 < q < \infty$,
\begin{align*}
\|K\|_q = \left(\int_{\mathbb{R}^d} |K(x)| \; |K(x)|^{q-1} dx \right)^{1/q} \leq \|K\|_\infty^{(q-1)/q} \|K\|_1^{1/q} <\infty.
\end{align*}

%
%
%
%

{\em Step 1.} Let $B_n(x) = \mathbb{E}\wh f(x) - f(x)$ and $\sigma_n(x)=\sqrt{\text{Var}(\wh f(x))}$. We will first prove the following facts, which show that $s_n^2$ and $\beta_{\mathbf{h}}(x)$ are the asymptotic expressions of the variance and bias of kernel density estimation uniformly in a small neighborhood of $\mathcal{M}$.
\begin{align}
&\sup_{x\in\mathcal{M}\oplus\epsilon_n} | B_n(x)-\beta_{\mathbf{h}}(x)| = o(\|\mathbf{h}\|^\nu),\label{fact2}\\
&\sup_{x\in\mathcal{M}\oplus\epsilon_n} |\sigma_n^2(x) - s_n^2| = o\left(\frac{1}{nh_1\cdots h_d}\right).\label{fact1}
\end{align}

We first show (\ref{fact2}). Note that by using the Taylor expansion for $f(x-\mathbf{h}\odot y)$ around $f(x)$, we have 
\begin{align*}
B_n(x) = & \int_{\mathbb{R}^d} [f(x-\mathbf{h}\odot y) - f(x)]K(y) dy \\
= & \int_{\mathbb{R}^d} [f(x+\mathbf{h}\odot y) - f(x)]K(y) dy \\
= & \sum_{\mathbf{i}\in\mathbb{Z}_+^{\nu,d}} \mathbf{h}^{(\mathbf{i})}  \int_{\mathbb{R}^d}\int_0^1 \frac{(1-t)^{\nu-1}}{(\nu-1)!} f_{(\mathbf{i})}(x+t\mathbf{h}\odot y) dt\, K(y)y^{(\mathbf{i})} dy.
%
\end{align*}
The assumption that $K$ is a $\nu$th order symmetric kernel function implies that we can write
\begin{align*}
\beta_{\mathbf{h}}(x) = \int_0^1 \frac{(1-t)^{\nu-1}}{(\nu-1)!}dt \sum_{\mathbf{i}\in\mathbb{Z}_+^{\nu,d}} \left[\mathbf{h}^{(\mathbf{i})} f_{(\mathbf{i})}(x)\, \int_{\mathbb{R}^d} K(y)y^{(\mathbf{i})} dy \right].
\end{align*}
Notice that 
\begin{align*}
& \sup_{x\in\mathcal{M}\oplus\epsilon_n}|B_n(x)- \beta_{\mathbf{h}}(x) | \nonumber\\
\leq & \sum_{\mathbf{i}\in\mathbb{Z}_+^{\nu,d}} \mathbf{h}^{(\mathbf{i})} \int_{\mathbb{R}^d}|K(y) y^{(\mathbf{i})} |\int_0^1 \frac{(1-t)^{\nu-1}}{(\nu-1)!} \sup_{x\in\mathcal{M}\oplus\epsilon_n} |f_{(\mathbf{i})}(x+t\mathbf{h}\odot y) - f_{(\mathbf{i})}(x)| dt dy.
%
\end{align*} 
We then obtain (\ref{fact2}) by applying the Dominated Convergence Theorem and assumption (F1) to the right-hand side of the above inequality.\\

Next we show (\ref{fact1}). For $s\in[-\epsilon_n,\epsilon_n]$ with the same $\epsilon_n$ given in (\ref{epsilon_n}), by using Taylor expansion, we have that 
\begin{align}\label{fact0}
\sup_{x\in\mathcal{M}} \left|f(\zeta_x(s)) - c - s\|\nabla f(x)\| \right| = O(s^2) = O(\epsilon_n^2).
\end{align}
Therefore using Taylor expansion again we have
\begin{align*}
& \sup_{x\in\mathcal{M}\oplus\epsilon_n} |\sigma_n^2(x) - s_n^2| \\
= &\sup_{x\in\mathcal{M}\oplus\epsilon_n} \left|\frac{1}{nh_1\cdots h_d} \int_{\mathbb{R}^d} [f(x-\mathbf{h}\odot y)-c] K^2(y) dy \right.\\
&\hspace{2cm} \left.- \frac{1}{n}\left[\int_{\mathbb{R}^d}  f(x-\mathbf{h}\odot y) K(y)dy\right]^2\right|\\
\leq & \frac{1}{nh_1\cdots h_d} \sup_{x\in\mathbb{R}^d} \|\nabla f(x) \|\left[ \epsilon_n + \|\mathbf{h}\| \int_{\mathbb{R}^d} \|y\|K^2(y)dy \right] + \frac{1}{n} \|f\|_\infty^2 \|K\|_1^2 \\ 
= & o\left(\frac{1}{nh_1\cdots h_d}\right).
\end{align*}
To get the above $o$-term, we have used the fact that 
\begin{align*}
\int_{\mathbb{R}^d} \|y\|K^2(y)dy = & \int_{\mathcal{B}_{\mathbf{0}}(1)} \|y\| K^2(y) dy + \int_{[\mathcal{B}_{\mathbf{0}}(1)]^\complement} \|y\| K^2(y) dy \\
\leq & \|K\|_2^2 + \|K\|_\infty \int_{\mathbb{R}^d} \|y\|^\nu |K(y)| dy \\
\leq & \|K\|_2^2 + \|K\|_\infty d^{\nu/2-1} \int_{\mathbb{R}^d}  (|y_1|^\nu +\cdots + |y_d|^\nu) |K(y)| dy <\infty,
\end{align*}
where $\mathbf{0}$ is the origin of $\mathbb{R}^d$, and we use (\ref{Kproperty}) and the definition of $\nu$th order kernels.

{\em Step 2.} We prove (\ref{L1expexpress}) in this step. Note that
\begin{align*}
\mathbb{E} \lambda_g(\mathcal{L} \Delta \widehat{\mathcal{L}}) & = \mathbb{E} \int \mathbf{I}(x\in\mathcal{L} \Delta \widehat{\mathcal{L}}) g(x)dx\\
& =  \int \mathbb{P}(x\in\mathcal{L} \Delta \widehat{\mathcal{L}}) g(x)dx\\
%
%
& = \int \mathbb{P}(\wh f(x) \geq c > f(x)) g(x) dx +  \int \mathbb{P}( f(x) \geq c > \wh f(x)) g(x)dx\\
& = \int_{\mathcal{L}^\complement} \mathbb{P}(\wh f(x) \geq c) g(x) dx + \int_{\mathcal{L}} \mathbb{P}(\wh f(x) < c) g(x) dx .
\end{align*}
Since $\mathcal{L}\Delta\widehat{\mathcal{L}} \subset \mathcal{M}\oplus\epsilon_n$ for large $n$ with probability one,
\begin{align}\label{symmexpress}
&\mathbb{E} \lambda_g(\mathcal{L} \Delta \widehat{\mathcal{L}}) \nonumber \\
= &\int_{\mathcal{L}^\complement \cap (\mathcal{M}\oplus\epsilon_n)} \mathbb{P}(\wh f(x) \geq c) g(x) dx + \int_{\mathcal{L}\cap (\mathcal{M}\oplus\epsilon_n)} \mathbb{P}(\wh f(x) < c) g(x) dx.
\end{align}
Note that here $\mathcal{L}^\complement \cap (\mathcal{M}\oplus\epsilon_n) = \{\zeta_{x}(s): \; x\in\mathcal{M}, -\epsilon_n<s<0\}$ and $\mathcal{L}\cap (\mathcal{M}\oplus\epsilon_n) = \{\zeta_{x}(s): \; x\in\mathcal{M}, 0\leq s<-\epsilon_n\}$. Similar to (\ref{reducedset2}), we have
\begin{align}\label{approxexpress1}
\mathbb{E} \lambda_g(\mathcal{L} \Delta \widehat{\mathcal{L}}) = & \left[\int_{\mathcal{M} } \int_{-\epsilon_n}^0 \mathbb{P} ( \wh f(\zeta_{x}(s))  \geq c ) g(\zeta_{x}(s)) ds d\mathscr{H}(x) \right. \nonumber\\
&\hspace{1cm} \left.+ \int_{\mathcal{M} } \int_0^{\epsilon_n} \mathbb{P} ( \wh f(\zeta_{x}(s))  < c ) g(\zeta_{x}(s)) ds d\mathscr{H}(x) \right] (1+o(1)).
\end{align}
For the leading term on the right-hand side of the above expression, it follows from (\ref{gapprox}) that
\begin{align}\label{approxexpress2}
& \int_{\mathcal{M} } \int_{-\epsilon_n}^0 \mathbb{P} ( \wh f(\zeta_{x}(s))  \geq c ) g(\zeta_{x}(s)) ds d\mathscr{H}(x) \nonumber\\
&\hspace{1cm} + \int_{\mathcal{M} } \int_0^{\epsilon_n} \mathbb{P} ( \wh f(\zeta_{x}(s))  < c ) g(\zeta_{x}(s)) ds d\mathscr{H}(x) \nonumber\\
= & \left[ \int_{\mathcal{M} } g^{(p)}(x) \int_{-\epsilon_n}^0 \mathbb{P} ( \wh f(\zeta_{x}(s))  \geq c ) |s|^p ds d\mathscr{H}(x) \right. \nonumber\\
&\hspace{1cm} \left.+ \int_{\mathcal{M} } g^{(p)}(x) \int_0^{\epsilon_n} \mathbb{P} ( \wh f(\zeta_{x}(s))  < c )  |s|^p ds d\mathscr{H}(x) \right] (1+o(1)).
\end{align}
Only focusing on the leading term again, in what follows we perform a sequence of decompositions. In general, we use $L_1$, $L_2$, and $L_3$ to denote dominant terms and $R_1$, $R_2$, and $R_3$ to denote remainder terms. We have
\begin{align}
& \int_{\mathcal{M} } g^{(p)}(x) \int_{-\epsilon_n}^0 \mathbb{P} ( \wh f(\zeta_{x}(s))  \geq c ) |s|^p ds d\mathscr{H}(x) \nonumber\\
&\hspace{1cm} + \int_{\mathcal{M} } g^{(p)}(x) \int_0^{\epsilon_n} \mathbb{P} ( \wh f(\zeta_{x}(s))  < c )  |s|^p ds d\mathscr{H}(x) \nonumber\\
& = \int_{\mathcal{M} } g^{(p)}(x) \int_{-\epsilon_n}^0 \mathbb{P} \left( \frac{\wh f(\zeta_{x}(s)) - \mathbb{E} \wh f(\zeta_{x}(s))}{\sigma_n(\zeta_{x}(s))} \geq \frac{c - \mathbb{E} \wh f(\zeta_{x}(s))}{\sigma_n(\zeta_{x}(s))} \right)  |s|^p ds d\mathscr{H}(x) \nonumber \\
&\hspace{0.5cm} + \int_{\mathcal{M} } g^{(p)}(x) \int_0^{\epsilon_n} \mathbb{P} \left( \frac{\wh f(\zeta_{x}(s)) - \mathbb{E} \wh f(\zeta_{x}(s))}{\sigma_n(\zeta_{x}(s))} < \frac{c - \mathbb{E} \wh f(\zeta_{x}(s))}{\sigma_n(\zeta_{x}(s))} \right)  |s|^p ds d\mathscr{H}(x)\nonumber\\
& =  L_1 + R_1,\label{R1term}
\end{align}
where 
\begin{align*}
& L_1 = \int_{\mathcal{M} } g^{(p)}(x)  \int_{-\epsilon_n}^0 \Phi \left( - \frac{c - \mathbb{E} \wh f(\zeta_{x}(s))}{\sigma_n(\zeta_{x}(s))} \right) |s|^p ds d\mathscr{H}(x) \\
&\hspace{1cm} + \int_{\mathcal{M} } g^{(p)}(x) \int_0^{\epsilon_n} \Phi \left(  \frac{c - \mathbb{E} \wh f(\zeta_{x}(s))}{\sigma_n(\zeta_{x}(s))} \right)  |s|^p ds d\mathscr{H}(x),
\end{align*}
with $\Phi$ the standard normal distribution function. We will show in {\em step 4} that
\begin{align}\label{R1rate}
|R_1| = o(\|\mathbf{h}\|^{\nu(p+1)} + s_n^{p+1}).
\end{align}

Using the results given in (\ref{fact1}) and (\ref{fact2}), we have
\begin{align}
%
%
L_1 = L_2 + R_2\label{R2term},
\end{align}
where 
\begin{align*}
& L_2 = \int_{\mathcal{M} } g^{(p)}(x) \int_{-\epsilon_n}^0 \Phi \left( \frac{ s\|\nabla f(x)\| + \beta_{\mathbf{h}}(x)}{s_n} \right) |s|^p ds d\mathscr{H}(x) \nonumber\\
&\hspace{1cm} + \int_{\mathcal{M} } g^{(p)}(x) \int_0^{\epsilon_n} \Phi \left(  \frac{ - s\|\nabla f(x)\| -\beta_{\mathbf{h}}(x)}{s_n} \right)  |s|^p ds d\mathscr{H}(x). 
\end{align*}
We will show in {\em step 4} that
\begin{align}\label{R2rate}
|R_2| = o( \|\mathbf{h}\|^{\nu(p+1)} + s_n^{p+1}).
\end{align}

Let $u = s/s_n$. Then we continue to decompose $L_2$ as follows.  
\begin{align}
%
L_2 = L_3 + R_3 \label{R3term},
\end{align}
where
\begin{align*}
& L_3 = s_n^{p+1} \int_{\mathcal{M} } g^{(p)}(x) \int_{-\infty}^0 \Phi \left( u\|\nabla f(x)\| + \frac{ \beta_{\mathbf{h}}(x)}{s_n} \right)  |u|^p du d\mathscr{H}(x) \nonumber\\
& \hspace{1cm} + s_n^{p+1} \int_{\mathcal{M} } g^{(p)}(x) \int_0^{\infty} \Phi \left(  - u \|\nabla f(x)\| - \frac{ \beta_{\mathbf{h}}(x)}{s_n}\right)  |u|^p du d\mathscr{H}(x) .
\end{align*}
We will show in {\em step 4} that
\begin{align}\label{R3rate}
|R_3| = o( s_n^{p+1}).
\end{align}
Using integration by parts we can write $L_3$ as
\begin{align*}
%
&(-1)^p s_n^{p+1} \int_{\mathcal{M} } g^{(p)}(x)\int_{-\infty}^0 \Phi \left( u\|\nabla f(x)\| + \frac{ \beta_{\mathbf{h}}(x)}{s_n} \right)  u^p du d\mathscr{H}(x)\\
&\hspace{1cm} + s_n^{p+1} \int_{\mathcal{M} } g^{(p)}(x)\int_0^{\infty} \Phi \left(  - u \|\nabla f(x)\| - \frac{ \beta_{\mathbf{h}}(x)}{s_n}\right)  u^p du d\mathscr{H}(x)\\
= &\frac{(-1)^{p+1}}{p+1}s_n^{p+1} \int_{\mathcal{M} } g^{(p)}(x) \|\nabla f(x)\| \int_{-\infty}^0 \phi \left( u\|\nabla f(x)\| + \frac{ \beta_{\mathbf{h}}(x)}{s_n}\right)  u^{p+1} du d\mathscr{H}(x)\\
&+ \frac{1}{p+1}s_n^{p+1}\int_{\mathcal{M} } g^{(p)}(x)\|\nabla f(x)\| \int_{-\infty}^0 \phi \left( u\|\nabla f(x)\| + \frac{ \beta_{\mathbf{h}}(x)}{s_n}\right)  u^{p+1} du d\mathscr{H}(x)\\
=&\frac{1}{p+1}s_n^{p+1} \int_{\mathcal{M} } g^{(p)}(x)  \|\nabla f(x)\| \int_{-\infty}^\infty \phi \left( u\|\nabla f(x)\| + \frac{ \beta_{\mathbf{h}}(x)}{s_n} \right) |u|^{p+1} du d\mathscr{H}(x), 
\end{align*}
where $\phi$ is the pdf of a standard normal distribution. Using the variable transformation $v=u\|\nabla f(x)\| + \beta_{\mathbf{h}}(x)/s_n$, then we have
\begin{align}\label{L3express}
L_3 
=&\frac{1}{p+1}s_n^{p+1} \int_{\mathcal{M} } \frac{g^{(p)}(x) }{\|\nabla f(x)\|^{p+1}} \int_{-\infty}^\infty \phi(v) \left| v- \frac{ \beta_{\mathbf{h}}(x)}{s_n} \right|^{p+1} dv d\mathscr{H}(x)\nonumber\\
=&\frac{1}{p+1} \int_{\mathcal{M} } \frac{g^{(p)}(x) }{\|\nabla f(x)\|^{p+1}} \int_{-\infty}^\infty \phi(v) \left| s_n \, v- \beta_{\mathbf{h}}(x) \right|^{p+1} dv d\mathscr{H}(x)\nonumber\\
=&\frac{1}{p+1} \mathbb{E} \int_{\mathcal{M} } \frac{g^{(p)}(x) }{\|\nabla f(x)\|^{p+1}} \left| s_n\,Z + \beta_{\mathbf{h}}(x) \right|^{p+1} d\mathscr{H}(x), 
\end{align}
where $Z$ is a standard normal random variable, and in the last equality above we have used the symmetry of $Z$'s distribution. Note that 
\begin{align}\label{normlowerbound}
\mathbb{E}\left| s_n\,Z + \beta_{\mathbf{h}}(x) \right|^{p+1} \leq 2^{p} s_n^{p+1} \mathbb{E}|Z|^{p+1} + 2^{p}|\beta_{\mathbf{h}}(x)|^{p+1}.
\end{align}
Since $\mathcal{M}$ is a compact set and $\sup_{x\in\mathcal{M}}|\beta_{\mathbf{h}}(x)| \leq C_0 \|\mathbf{h}\|^\nu$ for some constant $C_0>0$, we obtain that
\begin{align}\label{L3rate}
L_3 = O(s_n^{p+1} + \|\mathbf{h}\|^{\nu(p+1)} ).
\end{align}
Then (\ref{L1expexpress}) follows from (\ref{symmexpress}), (\ref{approxexpress1}), (\ref{approxexpress2}),  (\ref{R1term}), (\ref{R1rate}), (\ref{R2term}), (\ref{R2rate}), (\ref{R3term}), (\ref{R3rate}), (\ref{L3express}) and (\ref{L3rate}).

{\em Step 3.} Now we prove (\ref{L1express2}), which is implied by
\begin{align}\label{L1express2imp}
& \mathbb{E} \int_{\mathcal{M} } \frac{g^{(p)}(x) }{\|\nabla f(x)\|^{p+1}} \left| \wh f(x) -f(x) \right|^{p+1} d\mathscr{H}(x) \nonumber\\
= &\mathbb{E} \int_{\mathcal{M} } \frac{g^{(p)}(x) }{\|\nabla f(x)\|^{p+1}} \left| s_n\,Z + \beta_{\mathbf{h}}(x) \right|^{p+1} d\mathscr{H}(x) + o(s_n^{p+1} + \|\mathbf{h}\|^{\nu(p+1)}).
\end{align}
We will apply Lemma 1 in Horv\'{a}th (\cite{Horvath91}) (see Lemma~\ref{HorvathLemma} in the appendix) to show (\ref{L1express2imp}).
Let $Y_i(x)=(h_1\cdots h_d)^{-1} K(\mathbf{h}^{-1}\odot(x-X_i))$. Then $\wh f(x)=n^{-1}\sum_{i=1}^n Y_i(x)$ and $Var(Y_i(x)) = n\sigma_n^2(x)$. Now with $w=B_n(x)/\sigma_n(x)$, we have
\begin{align}\label{horvathexpress}
\mathbb{E}\left| \wh f(x) -f(x) \right|^{p+1} = \frac{1}{n^{p+1}} \mathbb{E}\left|\sum_{i=1}^n[Y_i(x) - \mathbb{E} Y_i(x)] + n^{1/2} [n^{1/2} \sigma_n(x)] w \right|^{p+1}.
\end{align}
For $1\leq k\leq p+3$ and $x\in\mathcal{M}$, using the substitution $u=\mathbf{h}^{-1}\odot(x-y)$ we can write
\begin{align*}
\mathbb{E}|Y_1(x)|^k & = (h_1\cdots h_d)^{-k} \int_{\mathbb{R}^d} |K(\mathbf{h}^{-1}\odot(x-y))|^k f(y)dy\\
& = (h_1\cdots h_d)^{-(k-1)} \int_{\mathbb{R}^d} |K(u)|^k f(x-\mathbf{h} \odot u)du\\
& = (h_1\cdots h_d)^{-(k-1)} c \|K\|_k^k\{1+o(1)\},
\end{align*}
where the last step is a consequence of (\ref{Kproperty}), assumption (F1), and the Dominated Convergence Theorem. Since $\mathbb{E}Y_1(x) = c+o(1)$ uniformly in $x\in\mathcal{M}$ (see (\ref{fact2})), and for $2\leq k \leq p+3$, 
\begin{align*}
\sum_{j=0}^k (-1)^{k-j}\binom{k}{j}\mathbb{E}|Y_1(x)|^j |\mathbb{E}Y_1(x)|^{k-j} & \leq \mathbb{E} |Y_1(x) - \mathbb{E}Y_1(x)|^k\\
& \leq \sum_{j=0}^k \binom{k}{j}\mathbb{E}|Y_1(x)|^j |\mathbb{E}Y_1(x)|^{k-j},
\end{align*}
we obtain that for $x\in\mathcal{M}$ and $2\leq k \leq p+3$, 
\begin{align}\label{moments}
\mathbb{E} |Y_1(x) - \mathbb{E}Y_1(x)|^k = (h_1\cdots h_d)^{-(k-1)} c \|K\|_k^k \{1+o(1)\},
\end{align}
where $o(1)$ is uniform in $x\in\mathcal{M}\oplus\epsilon_n$.

By applying  Lemma~\ref{HorvathLemma} and using (\ref{fact2}), (\ref{fact1}), (\ref{horvathexpress}) and (\ref{moments}), there exist positive constants $C_1$, $C_2$ and $C_3$ such that for all $x\in\mathcal{M}$,
\begin{align}
& \left| \mathbb{E} | \wh f(x) -f(x) |^{p+1} -  \mathbb{E}\left| \sigma_n(x)Z +  B_n(x) \right|^{p+1} \right| \nonumber\\
\leq & \frac{1}{n^{p+1}} C_1 \left( 1+\left|\frac{B_n(x)}{\sigma_n(x)}\right|^{p}\right) \nonumber \\
&\hspace{1cm} \times\left\{ n^{p/2} [n^{1/2}\sigma_n(x)]^{p-2} \frac{c\|K\|_3^3}{(h_1\cdots h_d)^2} + [n^{1/2}\sigma_n(x)]^{-2} \frac{c\|K\|_{p+3}^{p+3}}{(h_1\cdots h_d)^{p+2}} \right\} \nonumber\\
\leq & C_2 \left( 1+ \frac{(nh_1\cdots h_d)^{p/2}}{c^{p/2}\|K\|_2^p}  \left|\beta_{\mathbf{h}}(x)\right|^{p}\right) \nonumber\\
&\hspace{1cm} \times \left( \frac{c^{p/2} \|K\|_2^{p-2}\|K\|_3^3}{(nh_1\cdots h_d)^{p/2+1}} + \frac{\|K\|_{p+3}^{p+3}}{c\|K\|_2^4(nh_1\cdots h_d)^{p+1}}\right) \nonumber\\
\leq & C_3 (s_n^{p+2} + s_n^2|\beta_{\mathbf{h}}(x)|^p). \label{inequality1}
\end{align}

Let $\gamma_{n,\mathbf{h}} = s_n + \|\mathbf{h}\|^\nu$. Denote $A_n(x,Z)=\gamma_{n,\mathbf{h}}^{-1}[s_nZ + \beta_{\mathbf{h}}(x)]$ and $D_n(x,Z)=\gamma_{n,\mathbf{h}}^{-1}[(\sigma_n(x)-s_n)Z + B_n(x)-\beta_{\mathbf{h}}(x)]$. Then we can write
\begin{align}\label{adexpress}
&(\gamma_{n,\mathbf{h}})^{-(p+1)} \left| \mathbb{E}\left| \sigma_n(x)Z +  B_n(x)  \right|^{p+1} - \mathbb{E} \left| s_n\,Z + \beta_{\mathbf{h}}(x) \right|^{p+1} \right| \nonumber\\
= & \left| \mathbb{E}\left| A_n(x,Z) + D_n(x,Z)  \right|^{p+1} - \mathbb{E} \left| A_n(x,Z) \right|^{p+1} \right|.
\end{align}
Using the expressions of $s_n$ and $\beta_{\mathbf{h}}(x)$ in (\ref{snsq}) and (\ref{betahx}), we have that for $k=1,\cdots,2(p+1)$, $$\mathbb{E}|A_n(x,Z)|^k \leq 2^{k-1}\{|\gamma_{n,\mathbf{h}}^{-1} s_n|^k \mathbb{E}|Z|^k + |\gamma_{n,\mathbf{h}}^{-1} \beta_{\mathbf{h}}(x)|^k\}=O(1).$$ Similarly, we have that $\mathbb{E}|D_n(x,Z)|^k=o(1)$ for $k=1,\cdots,2(p+1)$, by using (\ref{fact2}) and (\ref{fact1}). Therefore, on the one hand,
\begin{align}\label{upperbound2}
&\mathbb{E}\left| A_n(x,Z) + D_n(x,Z)  \right|^{p+1} - \mathbb{E} \left| A_n(x,Z) \right|^{p+1} \nonumber\\
\leq & \sum_{j=0}^{p} \binom{p+1}{j} \mathbb{E}[|A_n(x,Z)|^j|D_n(x,Z)|^{p+1-j}] \nonumber\\
\leq &  \sum_{j=0}^{p} \binom{p+1}{j} \sqrt{\mathbb{E}[|A_n(x,Z)|^{2j}] \mathbb{E}[|D_n(x,Z)|^{2(p+1-j)}]} = o(1),
\end{align}
where we have used the Cauchy-Schwarz inequality. On the other hand,
\begin{align*}
&\left| A_n(x,Z) + D_n(x,Z)  \right|^{p+1} \\
\geq & \left| |A_n(x,Z)| - |D_n(x,Z)|  \right|^{p+1}\\
= &\begin{cases}
\sum_{j=0}^{p+1} \binom{p+1}{j} \|A_n(x,Z)|^j[-|D_n(x,Z)|]^{p+1-j} & \text{if } |A_n(x,Z)| \geq |D_n(x,Z)| \\
\sum_{j=0}^{p+1} \binom{p+1}{j} [-|A_n(x,Z)|]^j|D_n(x,Z)|^{p+1-j} & \text{if } |A_n(x,Z)| < |D_n(x,Z)|.
\end{cases}
\end{align*}
Therefore similar to (\ref{upperbound2}) we have
\begin{align}\label{lowerbound2}
& \mathbb{E}\left| A_n(x,Z) + D_n(x,Z)  \right|^{p+1} - \mathbb{E}\left| A_n(x,Z) \right|^{p+1} \nonumber\\
\geq & -2 \mathbb{E} |D_n(x,Z)| - \sum_{j=0}^{p} \binom{p+1}{j} \mathbb{E}\{|A_n(x,Z)|^j|D_n(x,Z)|^{p+1-j}\} = o(1).
\end{align}
It then follows from (\ref{upperbound2}) and (\ref{lowerbound2}) that the right-hand side of (\ref{adexpress}) is of order $o(1)$, which further implies that 
%
%
%
%
%
%
%
\begin{align}\label{inequality2}
\left| \mathbb{E}\left| \sigma_n(x)Z +  B_n(x)  \right|^{p+1} - \mathbb{E} \left| s_n\,Z + \beta_{\mathbf{h}}(x) \right|^{p+1} \right| = o(s_n^{p+1} + \|\mathbf{h}\|^{\nu(p+1)}).
\end{align}
Then (\ref{L1express2imp}) and hence (\ref{L1express2}) immediately follow from (\ref{inequality1}), (\ref{inequality2}) and the fact that the $\mathcal{M}$ is compact and $g^{(p)}(x)/\|\nabla f(x)\|^{p+1}$ is bounded on $\mathcal{M}$.

{\em Step 4}. We will prove (\ref{R1rate}), (\ref{R2rate}) and (\ref{R3rate}), as required in {\em step 2}.

We first show the proof of (\ref{R1rate}) for $R_1$. By the nonuniform Berry-Esseen theorem (c.f. Theorem 14, Petrov \cite{Petrov75}, page 125), there exists a constant $C_4>0$ such that for all $y\in\mathbb{R}$,
\begin{align*}
 & \left| \mathbb{P}\left( \frac{\wh f(\zeta_{x}(s)) - \mathbb{E}\wh f(\zeta_{x}(s)) }{\sigma_n(\zeta_{x}(s))} \leq y\right) - \Phi(y)\right| \\
 \leq & \frac{C_4 \mathbb{E} \left| Y_1(\zeta_{x}(s)) - \mathbb{E} Y_1(\zeta_{x}(s)) \right|^3}{n^{1/2}\left(\mathbb{E} \left| Y_1(\zeta_{x}(s)) - \mathbb{E} Y_1(\zeta_{x}(s)) \right|^2\right)^{3/2}(1+|y|)^3}.
\end{align*}
It then follows from (\ref{moments}) that there exists a constant $C_5>0$ such that for all $y\in\mathbb{R}$,
\begin{align}\label{fact3}
\sup_{x\in\mathcal{M}}\sup_{s\in[-\epsilon_n,\epsilon_n]} \left| \mathbb{P}\left( \frac{\wh f(\zeta_{x}(s)) - \mathbb{E}\wh f(\zeta_{x}(s)) }{\sigma_n(\zeta_{x}(s))} \leq y\right) - \Phi(y)\right| \leq \frac{C_5}{\sqrt{nh_1\cdots h_d}(1+|y|^3)}.
\end{align}
As a result,
\begin{align}\label{R1ineq}
|R_1| \leq \frac{C_5}{\sqrt{nh_1\cdots h_d}} \int_{\mathcal{M} } g^{(p)}(x)  \int_{-\epsilon_n}^{\epsilon_n} |s|^p\left(1+\left|  \frac{c - \mathbb{E} \wh f(\zeta_{x}(s))}{\sigma_n(\zeta_{x}(s))} \right|^3 \right)^{-1}   ds d\mathscr{H}(x).
\end{align}
Note that due to (\ref{fact2}), (\ref{fact1}) and (\ref{fact0}), there exists a positive constant $C_6$ such that for all $\eta_3\|\mathbf{h}\|^\nu \leq |s| \leq \epsilon_n$ (where $\eta_3$ appears in (\ref{epsilon_n})),
\begin{align*}
\inf_{x\in\mathcal{M}} \left|  \frac{c - \mathbb{E} \wh f(\zeta_{x}(s))}{\sigma_n(\zeta_{x}(s))} \right| \geq \frac{C_6|s|}{s_n}.
\end{align*}
Plugging this inequality to the right-hand side of (\ref{R1ineq}), we obtain
\begin{align*}
|R_1| \leq R_{11} + R_{12},
\end{align*}
where 
\begin{align}\label{R11rate}
R_{11} & = \frac{C_5}{\sqrt{nh_1\cdots h_d}} \int_{\mathcal{M} } g^{(p)}(x) d\mathscr{H}(x) \int_{|s| \leq \eta_3\|\mathbf{h}\|^\nu } |s|^p ds \nonumber \\
& = \frac{2C_5}{\sqrt{nh_1\cdots h_d}} \int_{\mathcal{M} } g^{(p)}(x) d\mathscr{H}(x) \frac{(\eta_3\|\mathbf{h}\|^\nu)^{p+1}}{p+1} \nonumber \\
& = o(\|\mathbf{h}\|^{\nu(p+1)}) ,
\end{align}
and 
\begin{align*}
R_{12} & =\frac{C_5}{\sqrt{nh_1\cdots h_d}} \int_{\mathcal{M} } g^{(p)}(x) d\mathscr{H}(x) \int_{\eta_3\|\mathbf{h}\|^\nu \leq |s| \leq \epsilon_n } |s|^p \left(1+\frac{C_6^3|s|^3}{s_n^3} \right)^{-1} ds .
\end{align*}
Using the variable transformation $t=s/s_n$ we have
\begin{align}\label{R12rate}
R_{12} & = \frac{2C_5s_n^{p+1}}{\sqrt{nh_1\cdots h_d}} \int_{\mathcal{M} } g^{(p)}(x) d\mathscr{H}(x) \int_{\eta_3\|\mathbf{h}\|^\nu/s_n \leq t \leq \epsilon_n/s_n} \frac{t^p}{1+C_6^3t^3}dt \nonumber\\
& \leq \frac{2C_5s_n^{p+1}}{\sqrt{nh_1\cdots h_d}} \int_{\mathcal{M} } g^{(p)}(x) d\mathscr{H}(x) \int_{0 \leq t \leq \epsilon_n/s_n} \frac{t^p}{1+C_6^3t^3}dt \nonumber\\
& = o(s_n^{p+1}+\|\mathbf{h}\|^{\nu(p+1)}),
\end{align}
where the last rate follows from assumption (H2) and 
\begin{align*}
\int_{0 \leq t \leq \epsilon_n/s_n} \frac{t^p}{1+C_6^3t^3}dt=
\begin{cases}
O(1), & p=0,1\\
O(\log(\epsilon_n/s_n)), & p=2 \\
O((\epsilon_n/s_n)^{p-2}), & p\geq3.
\end{cases}
\end{align*}
With (\ref{R11rate}) and (\ref{R12rate}), we thus get (\ref{R1rate}).

Next we show the proof of (\ref{R2rate}) for $R_2$. It follows from (\ref{fact2}), (\ref{fact1}) and (\ref{fact0}) that for any $\epsilon>0$ small enough, there exists $N_0>0$ such that for all $n>N_0$ we have that for all $|s|\leq\epsilon_n$
\begin{align*}
& \sup_{x\in\mathcal{M}} \left|  \frac{c - \mathbb{E} \wh f(\zeta_{x}(s))}{\sigma_n(\zeta_{x}(s))}  - \frac{ - s\|\nabla f(x)\| -\beta_{\mathbf{h}}(x)}{s_n}  \right| \\
\leq & \sup_{x\in\mathcal{M}} \left|  \frac{c - \mathbb{E} \wh f(\zeta_{x}(s))}{\sigma_n(\zeta_{x}(s))}  - \frac{c - \mathbb{E} \wh f(\zeta_{x}(s))}{s_n}   \right| +  \left|  \frac{c - \mathbb{E} \wh f(\zeta_{x}(s))}{s_n}   - \frac{ - s\|\nabla f(x)\| -\beta_{\mathbf{h}}(x)}{s_n}  \right|\\
\leq & \epsilon^2 \left(\frac{\|\mathbf{h}\|^\nu}{s_n} + \frac{|s|}{s_n}\right).
\end{align*}
Hence for large $n$, by possibly decreasing $\epsilon$ and increasing $\eta_3$ in (\ref{epsilon_n}) we have 
\begin{align*}
D_n^+(s,x) := &\left| \Phi \left( \frac{c - \mathbb{E} \wh f(\zeta_{x}(s))}{\sigma_n(\zeta_{x}(s))} \right) - \Phi \left( \frac{ - s\|\nabla f(x)\| -\beta_{\mathbf{h}}(x)}{s_n}  \right) \right| \\
\leq & 
\begin{cases}
1 & \text{ if } 0\leq s \leq \frac{\epsilon s_n - \beta_{\mathbf{h}}(x)}{\|\nabla f(x)\|} \\
 \epsilon^2 \left(\frac{\|\mathbf{h}\|^\nu}{s_n} + \frac{|s|}{s_n}\right) \phi \left( \frac{ s\|\nabla f(x)\| +\beta_{\mathbf{h}}(x) }{2s_n}  \right) & \text{ if } \frac{\epsilon s_n - \beta_{\mathbf{h}}(x)}{\|\nabla f(x)\|} < s \leq \epsilon_n
\end{cases},
\end{align*}
and similarly
\begin{align*}
D_n^-(s,x) := &\left| \Phi \left( - \frac{c - \mathbb{E} \wh f(\zeta_{x}(s))}{\sigma_n(\zeta_{x}(s))} \right) - \Phi \left( \frac{ s\|\nabla f(x)\| + \beta_{\mathbf{h}}(x)}{s_n}  \right) \right| \\
 \leq &
\begin{cases}
1 & \text{ if } 0\geq s \geq \frac{-\epsilon s_n - \beta_{\mathbf{h}}(x)}{\|\nabla f(x)\|} \\
 \epsilon^2 \left(\frac{\|\mathbf{h}\|^\nu}{s_n} + \frac{|s|}{s_n}\right) \phi \left( \frac{ s\|\nabla f(x)\| +\beta_{\mathbf{h}}(x) }{2s_n}  \right) & \text{ if } \frac{-\epsilon s_n - \beta_{\mathbf{h}}(x)}{\|\nabla f(x)\|} > s \geq -\epsilon_n
\end{cases}.
\end{align*}
Therefore
\begin{align}\label{R2decomp}
|R_2| \leq &\int_{\mathcal{M} } g^{(p)}(x) \left[ \int_{-\epsilon_n}^0 D_n^-(s,x) |s|^p ds + \int_0^{\epsilon_n} D_n^+(s,x)  |s|^p ds \right] d\mathscr{H}(x) \nonumber\\
\leq & R_{21} + R_{22},
\end{align}
where
\begin{align*}
R_{21} = \int_{\mathcal{M} } g^{(p)}(x) \int_{\frac{-\epsilon s_n - \beta_{\mathbf{h}}(x)}{\|\nabla f(x)\|}}^{\frac{\epsilon s_n - \beta_{\mathbf{h}}(x)}{\|\nabla f(x)\|}} |s|^p ds d\mathscr{H}(x),
\end{align*}
and
\begin{align*}
R_{22} = & \epsilon^2 \int_{\mathcal{M} } g^{(p)}(x) \int_{-\infty}^{\infty} \left(\frac{\|\mathbf{h}\|^\nu}{s_n} + \frac{|s|}{s_n}\right) \phi \left( \frac{ s\|\nabla f(x)\| +\beta_{\mathbf{h}}(x) }{2s_n}  \right) |s|^p ds d\mathscr{H}(x).
\end{align*}
%
%
Using the variable transformation $v =  [s\|\nabla f(x)\| + \beta_{\mathbf{h}}(x)]/{s_n}$, we get
\begin{align}
R_{21} = & s_n \int_{\mathcal{M} } g^{(p)}(x) \int_{-\epsilon}^\epsilon \frac{|s_nv-\beta_{\mathbf{h}}(x)|^p}{\|\nabla f(x)\|^{p+1}} dv d\mathscr{H}(x)\nonumber\\
\leq & \sum_{j=0}^p \binom{p}{j}s_n^{j+1}  \int_{\mathcal{M} } g^{(p)}(x) |\beta_{\mathbf{h}}(x)|^{p-j}  \int_{-\epsilon}^\epsilon \frac{|v|^j}{\|\nabla f(x)\|^{p+1}} dv d\mathscr{H}(x)\nonumber\\
= & \sum_{j=0}^p \frac{2}{j+1}  \binom{p}{j}s_n^{j+1} \epsilon^{j+1} \int_{\mathcal{M} }  \frac{g^{(p)}(x) |\beta_{\mathbf{h}}(x)|^{p-j} }{\|\nabla f(x)\|^{p+1}} d\mathscr{H}(x)\nonumber\\
\leq & \epsilon C_7 \left[ (\|\mathbf{h}\|^\nu)^{p+1} + s_n^{p+1}\right], \label{R2decomp1}
\end{align}
for some $C_7>0$. 
%
%
Using the variable transformation $v =  [s\|\nabla f(x)\| + \beta_{\mathbf{h}}(x)]/{s_n}$ again, we have
\begin{align}\label{R2decomp2}
R_{22} = &\epsilon^2 s_n \int_{\mathcal{M} } g^{(p)}(x)  \int_{-\infty}^{\infty} \left(\|\mathbf{h}\|^\nu + \frac{|s_nv-\beta_{\mathbf{h}}(x)|}{\|\nabla f(x)\|}\right) \phi \left( \frac{v}{2} \right)  \frac{|s_nv-\beta_{\mathbf{h}}(x)|^p}{\|\nabla f(x)\|^{p+1}} dv d\mathscr{H}(x) \nonumber\\
\leq & \epsilon^2 \sum_{j=0}^p \binom{p}{j}s_n^{j+1} \int_{\mathcal{M} } g^{(p)}(x) \frac{|\beta_{\mathbf{h}}(x)|^{p-j}}{\|\nabla f(x)\|^{p+1}} \int_{-\infty}^{\infty} \left(\|\mathbf{h}\|^\nu + \frac{s_n|v| + |\beta_{\mathbf{h}}(x)|}{\|\nabla f(x)\|}\right) \nonumber\\
& \hspace{1cm} \times\phi \left( \frac{v}{2} \right) |v|^j dv d\mathscr{H}(x) \nonumber\\
\leq &  \epsilon^2 C_8 \left[ (\|\mathbf{h}\|^\nu)^{p+2} + s_n^{p+2}\right],
\end{align}
for some $C_8>0$. 
%
%
%
Then (\ref{R2rate}) immediately follows from (\ref{R2decomp}), (\ref{R2decomp1}) and (\ref{R2decomp2}).

Next we show the proof of (\ref{R3rate}) for $R_3$. Note that $\epsilon_n/s_n\rightarrow\infty$ as $n\rightarrow\infty$ and $\sup_{x\in\mathcal{M}} |\beta_{\mathbf{h}}(x)| \leq C_9 \epsilon_n$ for some positive constant $C_9$. When $n$ is large enough,
\begin{align*}
|R_3| =& s_n^{p+1} \int_{\mathcal{M} } g^{(p)}(x) \int_{-\infty}^{-\epsilon_n/s_n} \Phi \left( u\|\nabla f(x)\| + \frac{ \beta_{\mathbf{h}}(x)}{s_n} \right)  |u|^p du d\mathscr{H}(x) \\
& \hspace{1cm} + s_n^{p+1} \int_{\mathcal{M} } g^{(p)}(x) \int_{\epsilon_n/s_n}^{\infty} \Phi \left(  - u \|\nabla f(x)\| - \frac{ \beta_{\mathbf{h}}(x)}{s_n}\right)  |u|^p du d\mathscr{H}(x)\\
\leq & 2 s_n^{p+1} \int_{\mathcal{M} } g^{(p)}(x) \int_{-\infty}^{-\epsilon_n/s_n} \Phi \left[ u(\|\nabla f(x)\| +C_9)\right]  |u|^p du d\mathscr{H}(x)\\
= & o(s_n^{p+1} ).
\end{align*}
Hence (\ref{R3rate}) is proved and here we conclude the proof. $\hfill\square$
\end{proof}

\begin{proof}[Proof of Proposition \ref{L2limit}] $\;$

An application of  Proposition A.1 in Cadre (\cite{Cadre06}) (see Lemma~\ref{CadreLemma} in the appendix) leads to
\begin{align*}
\int_{\mathcal{I}(\delta)} |\wh f(x) - f(x)|^2dx = \int_{c-\delta/2}^{c+\delta/2} \int_{f^{-1}(\tau)} \frac{|\wh f(x) - f(x)|^2}{\|\nabla f(x)\|} d\mathscr{H}(x)d\tau,
\end{align*}
for small $\delta>0$. Using the Lebesgue--Besicovitch theorem (cf. Evans and Gariepy \cite{Evans92}, Theorem 1, Chapter I), we obtain 
\begin{align}\label{L2AsLimit}
 \lim_{\delta\searrow0} \frac{1}{\delta}\mathbb{E} \int_{\mathcal{I}(\delta)} |\wh f(x) - f(x)|^2dx = \mathbb{E} \int_{\cal{M}} \frac{|\wh f(x) - f(x)|^2}{\|\nabla f(x)\|}d\mathscr{H}(x).
\end{align} 
Then the assertion follows from Theorem \ref{riskexpect}, where we take $p=1$ and $g^{(p)}(x)=c^r \|\nabla f(x)\|$ when $g(x)=f(x)^r|f(x)-c|$. See Remark~\ref{assumpremark} b)(ii). $\hfill\square$
\end{proof}

\begin{proof}[Proof of Theorem \ref{optimAMISE}] $\;$

Following simple algebra, we have for any $w>0$,
\begin{align*}
Q(\mathbf{u};\mathbf{M},a,\nu) = a^{2\nu/(d+2\nu)}w^{d/(d+2\nu)} Q(a^{-\nu/(d+2\nu)}w^{\nu/(d+2\nu)}\mathbf{u};w^{-1}\mathbf{M},1,\nu),
\end{align*}
and correspondingly,
\begin{align*}
\mathbf{u}(\mathbf{M},a,\nu) = a^{\nu/(d+2\nu)}w^{-\nu/(d+2\nu)} \mathbf{u}(w^{-1}\mathbf{M},1,\nu).
\end{align*}
The expression of the minimizer in (\ref{optimtheory}) then follows by noticing (\ref{AMISEexpress}) with $\mathbf{u} =\mathbf{h}^\nu$, $w=\kappa_\nu^2$, $\mathbf{M} = \kappa_\nu^2 A(f)$ and $a=n^{-1}cb(f)\|K\|_2^2$.

We continue to use the above notation in what follows. The argument of the uniqueness of the minimizer uses similar ideas in the proof of Theorem 6 in Yang and Tschernig (\cite{Yang99}), which we describe below. When $d=1$, for positive $\mathbf{u}$, 
\begin{align}\label{positiveconst}
\nabla^2 Q(\mathbf{u};\mathbf{M},a,\nu) =  \frac{2}{(\nu!)^2} \mathbf{M} + \frac{a(\nu+1)}{\nu^2} \mathbf{u}^{-2-1/\nu}>0.
\end{align}
When $d\geq2$, the Hessian of $Q(\mathbf{u};\mathbf{M},a,\nu)$ w.r.t. $\mathbf{u}$ is given by
\begin{align}\label{pdmatrix}
&\nabla^2 Q(\mathbf{u};\mathbf{M},a,\nu) \nonumber\\
= &\frac{2}{(\nu!)^2} \mathbf{M} + \frac{a}{\nu^2(u_1u_2\cdots u_d)^{1/\nu}} 
\begin{pmatrix} 
(\nu+1)u_1^{-2} & (u_1u_2)^{-1} &\cdots & (u_1u_d)^{-1} \\ 
(u_1u_2)^{-1}  & (\nu+1)u_2^{-2} &\cdots & (u_2u_d)^{-1} \\ 
\vdots & \vdots & \ddots & \vdots \\
(u_1u_d)^{-1} & (u_2u_d)^{-1} & \cdots & (\nu+1)u_d^{-2}
\end{pmatrix},
\end{align}
which is positive definite for $\mathbf{u}\in\bar{\mathbb{R}}_+^d$ under assumption (F2). Therefore by (\ref{AMISEexpress}), $\wt m(\mathbf{h})$ is a strictly convex function of $\mathbf{h}^\nu$ for $\mathbf{h}\in\bar{\mathbb{R}}_+^d$, which implies that there is at most one minimizer of $\wt m(\mathbf{h})$ in $\bar{\mathbb{R}}_+^d$. Also notice that $\wt m(\mathbf{h})$ tends to infinity if either $\|\mathbf{h}\| \rightarrow\infty$ or $\|\mathbf{h}\| \rightarrow0$. This shows that $\wt m(\mathbf{h})$ is uniquely minimized by $\wt{\mathbf{h}}_{\text{opt}}$. 

The result (\ref{hoptapprox}) follows a standard argument as given in Hall and Marron (\cite{Hall87}). We only sketch the proof here. By Fubini's theorem, we have from (\ref{MISEdef}) that
\begin{align*}
m({\bf h}) = \int_{\cal{M}} \frac{\mathbb{E} |\wh f(x) - f(x)|^2}{\|\nabla f(x)\|^2}d\mathscr{H}(x) = \int_{\cal{M}} \frac{\text{Var}(\wh f(x)) + |\mathbb{E} \wh f(x) - f(x)|^2}{\|\nabla f(x)\|^2}d\mathscr{H}(x).
\end{align*}
Under the assumption that $f$ has bounded and continuous $(\nu+2)$ times derivatives and $\int_{\mathbb{R}}|u^{\nu+2}\wt K(u)|du<\infty$, we can extend the expansions in (\ref{MISEbias}) and (\ref{MISEvar}) to the following:
\begin{align*}
&\mathbb{E} \wh f(x) -f(x) = \beta_{\mathbf{h}}(x)+ \frac{1}{(\nu+2)!}\kappa_{\nu+2}\sum_{k=1}^d h_k^{\nu+2} f_{(k * (\nu+2))}(x) + o(\|\mathbf{h}\|^{\nu+2}),\\
&\text{Var}(\wh f(x)) = s_n^2 + \frac{1}{nh_1\cdots h_d} \left( \frac{1}{2} \alpha(K) \sum_{k=1}^d h_k^2 f_{(k*2)}(x) +o(\|\mathbf{h}\|^2) \right) \\
&\hspace{2cm} - \frac{1}{n} \left( f(x) + \frac{1}{\nu!} k_\nu \sum_{k=1}^d h_k^\nu f_{(k*\nu)}(x) + o(\|\mathbf{h}\|^{\nu}) \right)^2,
\end{align*}
where $\kappa_{\nu+2} =\int_{\mathbb{R}} u^{\nu+2} \wt K(u) du$ and $\alpha(K) =\int_{\mathbb{R}^d} u_1^2 K(u)^2du<\infty$. Then we can obtain the following results: 
\begin{align}\label{mexpand}
&m(\mathbf{h}) = \wt m(\mathbf{h}) + O\left(\frac{1}{n\|\mathbf{h}\|^{d-2}} + \|\mathbf{h}\|^{2\nu+2}\right), \nonumber\\
& \nabla m(\mathbf{h})  = \nabla \wt m(\mathbf{h}) + O\left(\frac{1}{n\|\mathbf{h}\|^{d-1}} + \|\mathbf{h}\|^{2\nu+1}\right),\nonumber\\
& \nabla^2 m(\mathbf{h}) = \nabla^2 \wt m(\mathbf{h}) + O\left(\frac{1}{n\|\mathbf{h}\|^{d}} + \|\mathbf{h}\|^{2\nu}\right), \nonumber\\
\text{and } & \nabla^2 \wt m(\mathbf{h}) = O\left(\frac{1}{n\|\mathbf{h}\|^{d+2}} + \|\mathbf{h}\|^{2\nu-2}\right).
\end{align}
Using Taylor expansion, we have
\begin{align*}
0 = \nabla m(\mathbf{h}_{\text{opt}}) = \nabla m(\wt{\mathbf{h}}_{\text{opt}}) + \left[\int_0^1 \nabla^2 m(s\;\mathbf{h}_{\text{opt}} + (1-s)\;\wt{\mathbf{h}}_{\text{opt}}) ds \right](\mathbf{h}_{\text{opt}} - \wt{\mathbf{h}}_{\text{opt}}),
\end{align*}
which implies
\begin{align*}
\mathbf{h}_{\text{opt}} - \wt{\mathbf{h}}_{\text{opt}} = \left[\nabla^2\wt m(\wt{\mathbf{h}}_{\text{opt}}) \right]^{-1} \left[\nabla \wt m(\wt{\mathbf{h}}_{\text{opt}}) - \nabla m(\wt{\mathbf{h}}_{\text{opt}})\right] (1+o(1)),
\end{align*}
since $\nabla \wt m(\wt{\mathbf{h}}_{\text{opt}}) =0$ and $\nabla^2 \wt m(\wt{\mathbf{h}}_{\text{opt}})$ is a nonzero scalar when $d=1$ by (\ref{positiveconst}), or a nonsingular matrix when $d\geq2$ by (\ref{pdmatrix}). Then immediately we have (\ref{hoptapprox}) and (\ref{mhoptapprox}) by using (\ref{mexpand}).
$\hfill\square$
\end{proof}

\begin{proof}[Proof of Theorem \ref{practicald1}] $\;$

By noticing (\ref{stochrate}) and (\ref{biasuniform}), we have that with probability one $\wh N=N$ for $n$ sufficiently large (also see Theorem 3.1 in Biau et al., \cite{Biau07}). Hence we do not distinguish between $\wh N$ and $N$ in what follows. Without loss of generality, we assume that $x_1<\cdots < x_N$ and $\wh x_1<\cdots < \wh x_N$. Denote the kernel density estimation using bandwidth $h^{(k)}$ by $\wh f_k(x)$, $k=0,1, 2$. Also for $\ell=0,1,2,\cdots,$ denote the $\ell$th derivative of a function $g$ on $\mathbb{R}$ by $g^{(\ell)}$ if it exists, including the convention $g^{(0)}\equiv g$. 

Under assumption (F1), there exists $b_0>0$ such that $|f^\prime(x)|>\epsilon_0$ for $x\in\bigcup_{i=1}^N [x_i -b_0,x_i+b_0]$. By assuming $f$ has bounded continuous fourth derivatives and $K$ has bounded continuous third derivatives of bounded variation, it follows from Lemmas 2 and 3 in Arias-Castro et al. (\cite{Arias-Castro16}) that for $\wh f$ using bandwidth $h>0$ and for $\ell=0,1,2,$ and 3,
\begin{align}
&\sup_{x\in[x_i-b_0,x_i+b_0]}| \mathbb{E}\wh f^{(\ell)}(x) - f^{(\ell)}(x)| = O\left(h^{\min(4-\ell,2)} \right),\label{derivbiasorder}\\
&\sup_{x\in[x_i-b_0,x_i+b_0]}|\wh f^{(\ell)}(x) - \mathbb{E}\wh f^{(\ell)}(x)| = O\left(\sqrt{\frac{\log{n}}{nh^{1+2\ell}}} \right),\;a.s. \label{derivorder}
\end{align}
Due to the uniform consistency result for $\wh f_0$ shown in (\ref{hausdist}), (\ref{stochrate}) and (\ref{biasuniform}), with probability one $\wh x_i \in[x_i-b_0,x_i+b_0]$, for $i=1,\cdots,N$, for $n$ large enough. Let $g_n(x) = \wh f_0(x) - f(x)$. For $i=1,\cdots,N$, we have 
\begin{align}\label{gnxirate}
g_n(x_i)= O_p((nh^{(0)})^{-1/2} + (h^{(0)})^2) = O_p(n^{-2/5}).
\end{align}
Since $\wh f_0(\wh x_i) = f(x_i) = c$, we have 
\begin{align*}
g_n(\wh x_i) = \wh f_0(\wh x_i) - f(\wh x_i)= f(x_i) - f(\wh x_i) = f^\prime(\wt x_i) (x_i - \wh x_i),
\end{align*}
where $\wt x_i$ is between $x_i$ and $\wh x_i$ by using the Taylor expansion. Note that $\wt x_i \in [x_i-b_0, x_i + b_0]$, which implies that $|f^\prime(\wt x_i)|>\epsilon_0$ and yields
\begin{align}\label{fxdiff}
|x_i - \wh x_i| \leq \frac{1}{\epsilon_0} |g_n(\wh x_i)|. 
\end{align}
Another Taylor expansion for $g_n(\wh x_i)$ leads to
\begin{align}\label{fxdiff2}
g_n(\wh x_i) - g_n(x_i) = g_n^\prime (\check x_i)(\wh x_i - x_i),
\end{align}
where $\check x_i$ is between $\wh x_i$ and $x_i$ and 
\begin{align}\label{gnprimerate}
|g_n^\prime (\check x_i)| \leq \sup_{x\in[x_i-b_0,x_i+b_0]}|g_n^\prime(x)| = O_p\left(\sqrt{\frac{\log{n}}{n(h^{(0)})^3}} + (h^{(0)})^2\right) = o_p(1),
\end{align}
by using (\ref{derivbiasorder}) and (\ref{derivorder}). Then it follows from (\ref{gnxirate}), (\ref{fxdiff}), (\ref{fxdiff2}) and (\ref{gnprimerate}) that
\begin{align}
|\wh x_i - x_i|\leq \frac{1}{\epsilon_0}\left[ |g_n(x_i)| + |g_n(\wh x_i) - g_n(x_i)| \right]  = O_p(n^{-2/5}). 
\end{align}
Consequently, with $s_n(x): = \wh f_2^{\prime\prime}(x) - f^{\prime\prime}(x)$, using (\ref{derivbiasorder}) and (\ref{derivorder}) we have
\begin{align}
|s_n(\wh x_i) - s_n(x_i)| \leq |s_n^\prime(\breve x_i)|\; |(\wh x_i -x_i)|  = o_p(n^{-2/5}),
\end{align}
where $\breve x_i$ is between $\wh x_i$ and $x_i$.
With the choice $h^{(2)}=O(n^{-1/9})$,
\begin{align*} 
&\wh f_2^{\prime\prime}(\wh x_i) - f^{\prime\prime}(x_i) \\
= &[s_n(\wh x_i) - s_n(x_i)] + s_n(x_i)+ [f^{\prime\prime}(\wh x_i)- f^{\prime\prime}(x_i) ] = O_p(n^{-2/9}).
\end{align*}
Similarly, we have $\wh f_1^{\prime}(\wh x_i) - f^{\prime}(x_i) = O_p(n^{-2/7})$. We then have 
\begin{align*}
\sum_{i=1}^{N} [\wh f_2^{\prime\prime}(\wh x_i)]^2|\wh f_1^\prime(\wh x_i)|^{-1} - \sum_{i=1}^{N} [f^{\prime\prime}( x_i)]^2|f^\prime( x_i)|^{-1} = O_p(n^{-2/9}),
\end{align*}
and
\begin{align*}
 \sum_{i=1}^{N} |\wh f_1^\prime(\wh x_i)|^{-1} -  \sum_{i=1}^{N} |f^\prime(x_i)|^{-1} = O_p(n^{-2/7}).
\end{align*}
As a result, for $C$ and $\wh C$ given in (\ref{Constant}) and (\ref{optimprac1}), we have
\begin{align*}
\frac{\wh C}{C} -1 = O_p(h^{-2/9}).
\end{align*}
and correspondingly we get (\ref{hoptd1}) and (\ref{hoptd1m}). $\hfill\square$
\end{proof}

\section{Appendix}
In this appendix, we collect some known results that are used in the proofs. 

\begin{lemma}[Proposition A.1. in Cadre (\cite{Cadre06})]\label{CadreLemma}
Let $\phi: \mathbb{R}^d \mapsto \mathbb{R}_+$ be a continuously differentiable function such that $\phi(x)\rightarrow0$ as $\|x\|\rightarrow\infty$, and $J\subset\mathbb{R}_+$ be an interval such that $\inf J>0$ and $\inf_{\phi^{-1}(J)} \|\nabla \phi\|>0$. Then, for all bounded Borel function $g:\mathbb{R}^d\mapsto\mathbb{R}$:
\begin{align*}
\int_{\phi^{-1}(J)} gd\lambda =\int_J\int_{\phi^{-1}(\{s\})}\frac{g}{\|\nabla \phi\|}d\mathscr{H}ds.
\end{align*}
\end{lemma}

%
%
\begin{lemma}[Lemma 1 in Horv\'{a}th (\cite{Horvath91})]\label{HorvathLemma}
Let $Y,Y_1,\cdots,Y_n$ be i.i.d random vectors with $\mathbb{E}Y=\mu$, $var(Y)=\sigma^2$. If $\mathbb{E}(Y)^{p+2}<\infty$, then there is a constant $C=C(p)$ such that for any $w\in\mathbb{R}$,
\begin{align*}
& \left| \mathbb{E} \left| \sum_{i=1}^n(Y_i-\mu)+n^{1/2}\sigma w\right|^p - n^{p/2}\sigma^p\mathbb{E} |Z+w|^p\right| \\
\leq &C(1+|w|^{p-1}) \left[n^{(p-1)/2}\sigma^{p-3}\mathbb{E}|Y -\mu|^3 +\sigma^{-2}\mathbb{E}|Y-\mu|^{p+2} \right], 
\end{align*}
where $Z$ is a standard normal random variable.
\end{lemma}

\section*{Acknowledgement}
The author is grateful to Anand Vidyashankar, an Associate Editor and a referee for their careful reading of an earlier version of the paper and for insightful comments that lead to significant improvements. The research of Wanli Qiao was partially supported by NSF grants DMS 1821154 and FET 1900061, and a Jeffress Memorial Trust Award. The simulations in this work were run on ARGO, a research computing cluster provided by the Office of Research Computing at George Mason University, VA.



\end{document}